\documentclass[12pt]{article}
\usepackage{latexsym, amssymb}
\usepackage{graphicx}
\usepackage[rightcaption]{sidecap}
\usepackage{wrapfig}
\DeclareGraphicsExtensions{.png,.eps,.pdf}

\DeclareMathAlphabet\gothic{U}{euf}{m}{n}

\makeatletter
\def\eqnarray{\stepcounter{equation}\let\@currentlabel=\theequation
\global\@eqnswtrue
\tabskip\@centering\let\\=\@eqncr
$$\halign to \displaywidth\bgroup\hfil\global\@eqcnt\z@
  $\displaystyle\tabskip\z@{##}$&\global\@eqcnt\@ne
  \hfil$\displaystyle{{}##{}}$\hfil
  &\global\@eqcnt\tw@ $\displaystyle{##}$\hfil
  \tabskip\@centering&\llap{##}\tabskip\z@\cr}

\def\endeqnarray{\@@eqncr\egroup
      \global\advance\c@equation\m@ne$$\global\@ignoretrue}

\def\@yeqncr{\@ifnextchar [{\@xeqncr}{\@xeqncr[5pt]}}
\makeatother

\begin{document}
\bibliographystyle{tom}

\newtheorem{lemma}{Lemma}[section]
\newtheorem{thm}[lemma]{Theorem}
\newtheorem{cor}[lemma]{Corollary}
\newtheorem{voorb}[lemma]{Example}
\newtheorem{rem}[lemma]{Remark}
\newtheorem{prop}[lemma]{Proposition}
\newtheorem{stat}[lemma]{{\hspace{-5pt}}}
\newtheorem{obs}[lemma]{Observation}
\newtheorem{defin}[lemma]{Definition}

\newenvironment{remarkn}{\begin{rem} \rm}{\end{rem}}
\newenvironment{exam}{\begin{voorb} \rm}{\end{voorb}}
\newenvironment{defn}{\begin{defin} \rm}{\end{defin}}
\newenvironment{obsn}{\begin{obs} \rm}{\end{obs}}

\newenvironment{emphit}{\begin{itemize} }{\end{itemize}}

\newcommand{\gota}{\gothic{a}}
\newcommand{\gotb}{\gothic{b}}
\newcommand{\gotc}{\gothic{c}}
\newcommand{\gote}{\gothic{e}}
\newcommand{\gotf}{\gothic{f}}
\newcommand{\gotg}{\gothic{g}}
\newcommand{\gothh}{\gothic{h}}
\newcommand{\gotk}{\gothic{k}}
\newcommand{\gotm}{\gothic{m}}
\newcommand{\gotn}{\gothic{n}}
\newcommand{\gotp}{\gothic{p}}
\newcommand{\gotq}{\gothic{q}}
\newcommand{\gotr}{\gothic{r}}
\newcommand{\gots}{\gothic{s}}
\newcommand{\gotu}{\gothic{u}}
\newcommand{\gotv}{\gothic{v}}
\newcommand{\gotw}{\gothic{w}}
\newcommand{\gotz}{\gothic{z}}
\newcommand{\gotA}{\gothic{A}}
\newcommand{\gotB}{\gothic{B}}
\newcommand{\gotG}{\gothic{G}}
\newcommand{\gotL}{\gothic{L}}
\newcommand{\gotS}{\gothic{S}}
\newcommand{\gotT}{\gothic{T}}

\newcommand{\mn}{\marginpar{\hspace{1cm}*} }
\newcommand{\mnn}{\marginpar{\hspace{1cm}**} }

\newcommand{\mnq}{\marginpar{\hspace{1cm}*???} }
\newcommand{\mnnq}{\marginpar{\hspace{1cm}**???} }

\newcounter{teller}
\renewcommand{\theteller}{\Roman{teller}}
\newenvironment{tabel}{\begin{list}%
{\rm \bf \Roman{teller}.\hfill}{\usecounter{teller} \leftmargin=1.1cm
\labelwidth=1.1cm \labelsep=0cm \parsep=0cm}
                      }{\end{list}}

\newcounter{tellerr}
\renewcommand{\thetellerr}{(\roman{tellerr})}
\newenvironment{subtabel}{\begin{list}%
{\rm  (\roman{tellerr})\hfill}{\usecounter{tellerr} \leftmargin=1.1cm
\labelwidth=1.1cm \labelsep=0cm \parsep=0cm}
                         }{\end{list}}
\newenvironment{ssubtabel}{\begin{list}%
{\rm  (\roman{tellerr})\hfill}{\usecounter{tellerr} \leftmargin=1.1cm
\labelwidth=1.1cm \labelsep=0cm \parsep=0cm \topsep=1.5mm}
                         }{\end{list}}

\newcommand{\Ni}{{\bf N}}
\newcommand{\Ri}{{\bf R}}
\newcommand{\Ci}{{\bf C}}
\newcommand{\Si}{{\bf S}}
\newcommand{\Ti}{{\bf T}}
\newcommand{\Zi}{{\bf Z}}
\newcommand{\Fi}{{\bf F}}

\newcommand{\Epsilon}{{\rm E}}

\newcommand{\proof}{\mbox{\bf Proof} \hspace{5pt}} 
\newcommand{\remark}{\mbox{\bf Remark} \hspace{5pt}}
\newcommand{\ruimte}{\vskip10.0pt plus 4.0pt minus 6.0pt}

\newcommand{\simh}{{\stackrel{{\rm cap}}{\sim}}}
\newcommand{\ad}{{\mathop{\rm ad}}}
\newcommand{\Ad}{{\mathop{\rm Ad}}}
\newcommand{\Aut}{\mathop{\rm Aut}}
\newcommand{\arccot}{\mathop{\rm arccot}}
\newcommand{\capp}{{\mathop{\rm cap}}}
\newcommand{\rcapp}{{\mathop{\rm rcap}}}
\newcommand{\Capp}{{\mathop{\rm Cap}}}
\newcommand{\diam}{\mathop{\rm diam}}
\newcommand{\divv}{\mathop{\rm div}}
\newcommand{\dist}{\mathop{\rm dist}}
\newcommand{\codim}{\mathop{\rm codim}}
\newcommand{\RRe}{\mathop{\rm Re}}
\newcommand{\IIm}{\mathop{\rm Im}}
\newcommand{\Tr}{{\mathop{\rm Tr}}}
\newcommand{\Vol}{{\mathop{\rm Vol}}}
\newcommand{\card}{{\mathop{\rm card}}}
\newcommand{\supp}{\mathop{\rm supp}}
\newcommand{\sgn}{\mathop{\rm sgn}}
\newcommand{\essinf}{\mathop{\rm ess\,inf}}
\newcommand{\esssup}{\mathop{\rm ess\,sup}}
\newcommand{\Int}{\mathop{\rm Int}}
\newcommand{\Leibniz}{\mathop{\rm Leibniz}}
\newcommand{\lcm}{\mathop{\rm lcm}}
\newcommand{\loc}{{\rm loc}}

\newcommand{\mod}{\mathop{\rm mod}}
\newcommand{\spann}{\mathop{\rm span}}
\newcommand{\one}{1\hspace{-4.5pt}1}

\newcommand{\DWR}{}

\hyphenation{groups}
\hyphenation{unitary}

\newcommand{\tfrac}[2]{{\textstyle \frac{#1}{#2}}}

\newcommand{\cb}{{\cal B}}
\newcommand{\cc}{{\cal C}}
\newcommand{\cd}{{\cal D}}
\newcommand{\ce}{{\cal E}}
\newcommand{\cf}{{\cal F}}
\newcommand{\ch}{{\cal H}}
\newcommand{\ci}{{\cal I}}
\newcommand{\ck}{{\cal K}}
\newcommand{\cl}{{\cal L}}
\newcommand{\cm}{{\cal M}}
\newcommand{\cn}{{\cal N}}
\newcommand{\co}{{\cal O}}
\newcommand{\cs}{{\cal S}}
\newcommand{\ct}{{\cal T}}
\newcommand{\cx}{{\cal X}}
\newcommand{\cy}{{\cal Y}}
\newcommand{\cz}{{\cal Z}}

\newcommand{\wtozp}{W^{1,2}\raisebox{10pt}[0pt][0pt]{\makebox[0pt]{\hspace{-34pt}$\scriptstyle\circ$}}}
\newlength{\hightcharacter}
\newlength{\widthcharacter}
\newcommand{\covsup}[1]{\settowidth{\widthcharacter}{$#1$}\addtolength{\widthcharacter}{-0.15em}\settoheight{\hightcharacter}{$#1$}\addtolength{\hightcharacter}{0.1ex}#1\raisebox{\hightcharacter}[0pt][0pt]{\makebox[0pt]{\hspace{-\widthcharacter}$\scriptstyle\circ$}}}
\newcommand{\cov}[1]{\settowidth{\widthcharacter}{$#1$}\addtolength{\widthcharacter}{-0.15em}\settoheight{\hightcharacter}{$#1$}\addtolength{\hightcharacter}{0.1ex}#1\raisebox{\hightcharacter}{\makebox[0pt]{\hspace{-\widthcharacter}$\scriptstyle\circ$}}}
\newcommand{\scov}[1]{\settowidth{\widthcharacter}{$#1$}\addtolength{\widthcharacter}{-0.15em}\settoheight{\hightcharacter}{$#1$}\addtolength{\hightcharacter}{0.1ex}#1\raisebox{0.7\hightcharacter}{\makebox[0pt]{\hspace{-\widthcharacter}$\scriptstyle\circ$}}}

 \thispagestyle{empty}
  
 \begin{center}
\vspace*{1.5cm}

{\Large{\bf The weighted  Hardy inequality  }}\\[3mm]
{\Large{\bf and  self-adjointness  of }}\\[3mm]
{\Large{\bf symmetric diffusion operators }}  \\[5mm]
\large Derek W. Robinson$^\dag$ \\[1mm]

\normalsize{25th June 2020}\\[1mm]
\end{center}

\vspace{+5mm}

\begin{list}{}{\leftmargin=1.7cm \rightmargin=1.7cm \listparindent=15mm 
   \parsep=0pt}
   \item
{\bf Abstract} $\;$ 
Let $\Omega$ be a domain in $\Ri^d$ with  boundary $\Gamma$${\!,}$ $d_\Gamma$ the Euclidean distance to the boundary
 and $H=-\divv(C\,\nabla)$  an elliptic operator with
$C=(\,c_{kl}\,)>0$ where $c_{kl}=c_{lk}$ are  real, bounded, Lipschitz functions.
We assume that $C\sim c\,d_\Gamma^{\,\delta}$ as $d_\Gamma\to0$ in the sense of asymptotic analysis
where $c$ is a strictly positive, bounded, Lipschitz function and $\delta\geq0$.
We also assume that  there is an $r>0$ and a $ b_{\delta,r}>0$ such that the weighted Hardy inequality
\[
\int_{\Gamma_{\!\!r}} d_\Gamma^{\,\delta}\,|\nabla \psi|^2\geq b_{\delta,r}^{\,2}\int_{\Gamma_{\!\!r}} d_\Gamma^{\,\delta-2}\,| \psi|^2
\]
is valid for all $\psi\in C_c^\infty(\Gamma_{\!\!r})$ where $\Gamma_{\!\!r}=\{x\in\Omega: d_\Gamma(x)<r\}$.
We then prove that the condition $(2-\delta)/2<b_\delta$ is sufficient for the essential self-adjointness of $H$ on $C_c^\infty(\Omega)$ 
with  $b_\delta$  the supremum over $r$ of all possible $b_{\delta,r}$ in the Hardy inequality.
This result extends all known results for domains with smooth boundaries and also gives information on 
self-adjointness for a large family of domains with rough, e.g.\ fractal, boundaries.

\end{list}

\vfill

\noindent AMS Subject Classification: 31C25, 47D07.

\noindent Keywords: Self-adjointness, diffusion operators, weighted Hardy inequality.

\vspace{0.5cm}

\noindent
\begin{tabular}{@{}cl@{\hspace{10mm}}cl}
$ {}^\dag\hspace{-5mm}$&   Mathematical Sciences Institute (CMA)    &  {} &{}\\
  &Australian National University& & {}\\
&Canberra, ACT 0200 && {} \\
  & Australia && {} \\
  &derek.robinson@anu.edu.au
 & &{}\\
\end{tabular}

\newpage

\setcounter{page}{1}

\section{Introduction}\label{S1}

In this paper we derive sufficiency criteria for the self-adjointness of  degenerate elliptic operators 
$H=-\divv(C\nabla)$ defined on $C_c^\infty(\Omega)$  where $\Omega$ is a domain in $\Ri^d$ with boundary $\Gamma$
and  $C=(c_{kl})$ is a  strictly positive, symmetric, $d\times d$-matrix with  $c_{kl}$ real Lipschitz continuous functions.
We assume that $C$ resembles the diagonal matrix $c\,d_\Gamma^{\,\delta}I$
near the boundary where $c$ is  a strictly positive bounded Lipschitz function,  $d_\Gamma$  the Euclidean distance to the boundary 
and $\delta\geq0$ a parameter which measures the order of degeneracy.
A precise definition of the degeneracy condition will be given in Section~\ref{S5} but the idea is that the matrices
 $C$ and $c\,d_\Gamma^{\,\delta}I$ are equivalent on a boundary layer
 $\Gamma_{\!\!r}=\{x\in\Omega: d_\Gamma(x)<r\}$ and asymptotically equal as $r\to0$.
  In an earlier article \cite{Rob15} we  established that if $\Omega$ is a $C^2$-domain, 
  or the complement of a lower dimensional $C^2$-domain, 
  then the condition  $\delta>2-(d-d_{\!H})/2$, with  $d_{\!H}$  the Hausdorff dimension of the boundary $\Gamma$, is sufficient for self-adjointness.
 Broadly similar conclusions had been reached earlier for bounded domains  by Nenciu and Nenciu \cite{NeN} by quite disparate arguments.
 In this paper we develop an alternative approach  which reconciles the two sets of arguments
 and gives  sufficiency   criteria  for a much broader class of domains.
 In particular we are able to derive results for  domains with  rough boundaries, e.g.\  boundaries with a fractal 
 nature or boundaries which are uniformly disconnected.

  Our arguments rely on two basic ideas.
First,   self-adjointness is determined by the  properties of the coefficients of the operator on an arbitrarily thin boundary layer $\Gamma_{\!\!r}$.
  Secondly, the existence of the weighted Hardy inequality on the boundary layer,
  \begin{equation}
  \int_{\Gamma_{\!\!r}} d_\Gamma^{\,\delta}\,|\nabla \psi|^2\geq b_{\delta,r}^{\,2}\int_{\Gamma_{\!\!r}} d_\Gamma^{\,\delta-2}\,| \psi|^2
\label{esa1.2}
 \end{equation}
 for  some $b_{\delta,r}>0$ and all $\psi\in C_c^\infty(\Gamma_{\!\!r})$ is crucial.
 If the Hardy  inequality is valid for one $r>0$ it is clearly valid for all $s\in\langle0,r]$
 and one can choose the corresponding $b_{\delta,s}\geq b_{\delta,r}$.
Therefore  we define  the (boundary) Hardy constant $b_\delta$ by
 \[
 b_\delta= \textstyle{\sup_{r>0}}\;\textstyle{\inf_{\psi\in C_c^\infty(\Gamma_{\!\!r})}}\;\Big\{\Big( {\displaystyle\int_{\Gamma_{\!\!r}}} d_\Gamma^{\,\delta}\,|\nabla \psi|^2\Big)^{1/2}\Big
 /\,\Big({\displaystyle\int_{\Gamma_{\!\!r}}} d_\Gamma^{\,\delta-2}\,| \psi|^2\Big)^{1/2}\Big\}\;,
 \]
i.e.\  $b_\delta$ is the supremum over the possible choices of the $b_{\delta,s}$.
Then, under these two assumptions, 
 we establish in  Theorem~\ref{tsa5.1} that the condition
 \begin{equation}
 (2-\delta)/2<b_\delta
 \label{esa1.3}
 \end{equation}
is sufficient for self-adjointness of $H$.
It is notable that there are no explicit restrictions on the domain $\Omega$ or its boundary $\Gamma$ only implicit conditions necessary for  the boundary Hardy inequality (\ref{esa1.2}).
 In all the  specific cases considered in \cite{Rob15} the Hardy inequality~(\ref{esa1.2}) is valid  and  the Hardy constant has  the standard value $b_\delta=(d-d_{\!H}+\delta-2)/2$.
Then (\ref{esa1.3}) reduces to the  condition $\delta>2-(d-d_{\!H})/2$.
Consequently a similar conclusion is valid for all domains which support the  weighted Hardy inequality (\ref{esa1.2})  on a boundary layer with the standard constant.

Note that if $\delta\geq 2$  then (\ref{esa1.3}) is obviously satisfied. 
In fact in this case  self-adjointness of $H$ follows from an upper bound $C\leq a\,d_\Gamma^{\,\delta}I$ on  the coefficients
and the weighted Hardy inequality is irrelevant.
 This latter result was already known (see, for example, \cite{ERS5} Theorem~4.10 and Corollary~4.12)  but we give a short proof in Section~\ref{S2}.
 It is also known from an earlier collaboration \cite{LR} with Lehrb{\"a}ck on Markov uniqueness  that the condition $\delta\geq 2-(d-d_{\!H})$
 is necessary for self-adjointness.
Theorem~\ref{tsa5.1} establishes, however, that the condition $\delta>2-(d-d_{\!H})/2$ is sufficient for $H$ to be self-adjoint in the standard case.
The proof of this statement   utilizes the ideas of Agmon, \cite{Agm1} Theorem~1.5, as developed by Nenciu and Nenciu \cite{NeN}  but extended to unbounded domains.
These ideas are elaborated in detail in Section~\ref{S4} where we establish  a prototype of our  main result from the upper bounds $C\leq a\,d_\Gamma^{\,\delta}I$ 
and a stronger version of the weighted Hardy inequality (\ref{esa1.2}).
Then in Section~\ref{S5} we establish the key result,  Theorem~\ref{tsa5.1}.

In Section~\ref{S6} we apply the latter result to  uniform domains with Ahlfors regular boundaries.
This allows a wide range of `rough' boundaries.
The application is made possible by a result of Lehrb{\"a}ck,  \cite{Leh3} Theorem~1.3, which establishes the validity of the weighted Hardy inequality (\ref{esa1.2})
on unbounded John domains.
A modification of Lehrb{\"a}ck's arguments \cite{Leh5} also demonstrates the validity of the Hardy inequality on boundary layers for bounded John domains.
In combination with the Ahlfors regularity one can then deduce that the condition (\ref{esa1.3}) again suffices for the self-adjointness of $H$.
Unfortunately little is known about the optimal value of the (boundary)  Hardy constant $b_\delta$ at this level of generality.
Nevertheless we establish that $b_\delta$ is bounded above  by 
the standard value 
$(d-d_{\!H}+\delta-2)/2$ and that $b_\delta+\delta/2$ is an increasing function
on the interval of interest.
It then follows that there is a critical degeneracy $\delta_c\in\langle 2-(d-d_{\!H})/2, 2\rangle$ such that 
$H$ is self-adjoint if $\delta>\delta_c$.
In addition $b_\delta$  is equal to the standard value if and only if $b_2=(d-d_{\!H})/2$ and in this case $\delta_c=2-(d-d_{\!H})/2$
(see Theorem~\ref{tsa6.1}).
One can, however, construct examples for which $b_\delta<(d-d_{\!H}+\delta-2)/2$.
In fact for each $\delta>2-(d-d_{\!H})/2$ there are examples with
 $b_\delta$ arbitrarily close to zero  and consequently $\delta_c$
is arbitrarily close to $2$.
Moreover,   if the sufficiency condition (\ref{esa1.3}) is satisfied then $H$ satisfies a  weighted Rellich inequality on a boundary layer $\Gamma_{\!\!r}$.

Finally we  emphasize  that the weighted  Hardy inequality (\ref{esa1.2}) only depends on the operator  $H$ through the order of degeneracy $\delta$
and the conclusions of  our  theorems only depend  on properties near the boundary.
Thus  if $\Gamma$ decomposes as the countable union of positively separated components  $\Gamma^{(j)}$ Theorems~\ref{tsa5.1} and \ref{tsa6.1}   can be elaborated.
 In this situation the boundary layer $\Gamma_{\!\!r}$  also decomposes into separate components  $\Gamma^{(j)}_{\!\!r}\!\!,\,$ if $r$ is sufficiently small.
 Then one can assign different orders of degeneracy $\delta_j$  to each component and introduce different Hardy constants $b_{\delta_j}$.
After this modification  self-adjointness of $H$ follows   from the family of conditions $(2-\delta_j)/2<b_{\delta_j}$.
This is discussed in Section~4.

\section{Preliminaries}\label{S2}

In this section we gather some preliminary information on self-adjointness of diffusion operators.
The elliptic operator $H=-\divv(C\nabla)$,  with domain $C_c^\infty(\Omega)$,   is a positive symmetric operator on $L_2(\Omega)$.
Consequently it is closable with respect to the graph norm $\|\varphi\|_{D(H)}=(\|H\varphi\|_2^2+\|\varphi\|_2^2)^{1/2}$.
For simplicity of notation we let $H$ and $D(H)$ denote the closure of the operator and its domain, respectively.

Since the coefficients $c_{kl}$ are bounded $\sup_{x\in\Omega}\|C(x)\|<\infty$, i.e.\ there is a $\nu>0$ such that 
$C\leq  \nu I$.
It also follows from the strict positivity of $C$ that for each compact subset $K$ of $\Omega$ there is a $\mu_K>0$ such that $C\geq\mu_K I$.
Thus the operator $H$ is locally strongly elliptic.
These local properties  imply, by elliptic regularity, that the domain $D(H^*)$ of the $L_2$-adjoint $H^*$ of $H$ is contained in $W^{2,2}_{\rm loc}(\Omega)$.
Moreover, $W^{2,\infty}_{c}(\Omega)D(H^*)\subseteq D(H)$.

Next let $h$ denote the positive bilinear form associated with $H$ on $L_2(\Omega)$, i.e.\ the form with domain $D(h)=C_c^\infty(\Omega)$ given by 
\[
h(\psi,\varphi)=(\psi,H\varphi)=\sum\nolimits^d_{k,l=1}(\partial_k\psi, c_{kl}\partial_l\varphi)
\]
for all $\psi, \varphi\in D(h)$ and set $h(\varphi)=h(\varphi, \varphi)$.
The form is closable  with respect to the graph norm $\|\varphi\|_{D(h)}=(h(\varphi)+\|\varphi\|_2^2)^{1/2}$ and we  also use $h$ and $D(h)$ to denote the closed form and its domain.
Then $h$ is a Dirichlet form \cite{BH} \cite{FOT} and 
 by elliptic regularity $D(h)\subseteq W^{1,2}_{\rm loc}(\Omega)$.
The Dirichlet form $h$ has a {\it carr\'e du champ}, a positive bilinear form $\psi, \varphi\in D(h)\mapsto \Gamma_{\!\!c}(\psi, \varphi)\in L_1(\Omega)$ such that
\[
\Gamma_{\!\!c}(\psi, \varphi)=\sum\nolimits^d_{k,l=1}c_{kl} (\partial_k\psi)( \partial_l\varphi)
\]
for all $\psi, \varphi\in C_c^\infty(\Omega)$ (see, for example, \cite{BH} Section~1.4).
Consequently $h(\varphi)=\|\Gamma_{\!\!c}(\varphi)\|_1$ where $\Gamma_{\!\!c}(\varphi)=\Gamma_{\!\!c}(\varphi, \varphi)$.

\smallskip
In the remainder of this section we introduce a  well-known identity which has been used extensively  to establish  self-adjointness  and apply it to operators on $\Ri^d$ and to operators with $\delta\geq2$.
These applications are subsequently  useful.
At this stage we do not require any additional  restrictions on the boundary properties of the coefficients.
Further assumptions   will be introduced in  Section~\ref{S5}.

The standard Stone-von Neumann criterion for the self-adjointness of $H$ is the range property $R(\lambda I+H)=L_2(\Omega)$ for some, or for all,  $\lambda>0$.
Equivalently one has the kernel condition $\ker(\lambda I+H^*)=\{0\}$ on the adjoint.
The following fundamental proposition gives a method for verifying the latter condition in some very general situations.

\begin{prop}\label{psa2.1}
If $\varphi\in D(H^*)$ and $\eta\in W^{1,\infty}_c(\Omega)$ then $\eta\varphi\in D(h)$ and
\begin{equation}
(H^*\varphi, \eta^2\varphi)=h(\eta\varphi)-(\varphi,\Gamma_{\!\!c}(\eta)\varphi)
\;.\label{esa2.1}
\end{equation}
Thus if $(\lambda I+H^*)\varphi=0$ for some $\lambda>0$ then 
\begin{equation}
\lambda\, \|\eta\varphi\|_2^2+h(\eta\varphi)=(\varphi,\Gamma_{\!\!c}(\eta)\varphi)
\label{esa2.2}
\end{equation}
for all $\eta\in W^{1,\infty}_c(\Omega)$.
\end{prop}

This result has a long history and a variety of different proofs. 
The identity (\ref{esa2.1}) occurs in Wienholtz' 1958 thesis (see \cite{Wie} Section~3) so it was certainly known to Rellich in the 1950s.
It was also derived by Agmon in his 1982 lectures (see \cite{Agm1} equation (1.16)).
Both these authors used standard methods of elliptic differential equations.
Alternatively~(\ref{esa2.1}) can be established in the abstract setting of local Dirichlet forms, \cite{Rob12} Lemma~2.2,
or on graphs, \cite{KPP} Lemma~2.1.
The identity is referred to as the localization lemma in \cite{NeN} where numerous other background references are given.

The most straightforward application of the result follows by noting that the form $h$ is positive.
Hence if $(\lambda I+H^*)\varphi=0$ then (\ref{esa2.2}) implies that 
\begin{equation}
\lambda \,\|\eta\varphi\|_2^2\leq (\varphi,\Gamma_{\!\!c}(\eta)\varphi)
\label{esa2.3}
\end{equation}
for all $\eta\in W^{1,\infty}_c(\Omega)$.
This latter estimate can be exploited to deduce self-adjointness of $H$ by the construction of a sequence of $\eta$ which converges to the identity but
$\Gamma_{\!\!c}(\eta)$ converges to zero thereby implying that $\varphi=0$.
This argument appears as Theorem~3.1 of \cite{Dav14} who also gives several earlier references dating back to the 1960s.

\begin{remarkn}\label{rsa2.1}
There is a certain delicacy in the derivation of the identity (\ref{esa2.1})  since $D(H^*)$ is not generally a subset of $D(h)$.
Therefore it is essential that $\eta$ is a differentiable function with compact support.
This guarantees that $\eta D(H^*)\subseteq D(h)$ by elliptic regularity.
In fact $D(H^*)\subseteq D(h)$ if and only if $H$ is self-adjoint.
This is a consequence of a key property  of the Friedrichs extension $H_{\!F}$ of $H$.
This is the self-adjoint extension of $H$ determined by the Dirichlet form $h$ and it is the only self-adjoint extension with domain contained in $D(h)$ (see \cite{Kat1}
Theorem~VI.2.11).
Thus if $H$ is self-adjoint then $H^*=H=H_{\!F}$.
Consequently $D(H^*)\subseteq D(h)$.
Conversely all self-adjoint extensions of $H$ are restrictions of $H^*$ so if $D(H^*)\subseteq D(h)$ then they must all be equal to $H_{\!F}$ by the result cited by Kato.
\end{remarkn}

It is significant that the inequality (\ref{esa2.3}) does not explicitly depend on the form $h$.
Therefore the compactness of the support of $\eta$ is no longer critical.
Consequently the inequality can be extended by continuity to a larger class of functions.

\begin{cor}\label{csa2.1}
If $\varphi\in D(H^*)$ with  $(\lambda I+H^*)\varphi=0$ for some $\lambda>0$ then 
\begin{equation}
\lambda\, \|\eta\varphi\|_2^2\leq (\varphi,\Gamma_{\!\!c}(\eta)\varphi)
\label{esa2.4}
\end{equation}
for all $\eta\in\bigcap_{s>0} W^{1,\infty}(\Omega_s)$ where $\Omega_s=\{x\in\Omega: d_\Gamma(x)>s\}$.
\end{cor}
\proof\
Fix $\eta\in W^{1,\infty}(\Omega_s)$.
Then let $\rho\in W^{1,\infty}_c(\Ri^d)$ be a positive function with $\rho=1$ on a ball $B\subset \Ri^d$ and zero on the complement of a larger concentric ball.
Then define $\rho_n(x)=\rho(x/n)$.
Since $W^{1,\infty}_c(\Ri^d)$ is an algebra of multipliers on $W^{1,\infty}(\Omega)$ it follows that $\rho_n\eta\in  W^{1,\infty}_c(\Omega)$.
Therefore replacing $\eta$ in (\ref{esa2.2}) by $\rho_n\eta$ and using the Leibniz rule combined with the Cauchy--Schwarz inequality one deduces that for each 
$\varepsilon>0$ one has 
\[
\lambda\, \|\rho_n\eta\varphi\|_2^2\leq (\varphi,\Gamma_{\!\!c}(\rho_n\eta)\varphi)
\leq (1+\varepsilon)\,(\rho_n\varphi,\Gamma_{\!\!c}(\eta)\rho_n\varphi)+(1+\varepsilon^{-1})\,(\eta\varphi,\Gamma_{\!\!c}(\rho_n)\eta\varphi)
\;.
\]
But $ \|\rho_n\eta\varphi\|_2^2\to  \|\eta\varphi\|_2^2$ as $n\to\infty$.
Moreover,
\[
(\eta\varphi,\Gamma_{\!\!c}(\rho_n)\eta\varphi)\leq \nu\,(\eta\varphi,|\nabla\!\rho_n|^2\eta\varphi)\leq \nu\,\|\eta\varphi\|_2^2\,\|\nabla\!\rho\|_\infty^2/n^2\to0
\]
as $n\to\infty$ by the definition of $\rho_n$.
The conclusion follows immediately.
\hfill$\Box$

\bigskip

This argument will be used in a more complicated context in Section~\ref{S4}.
Another simple corollary of the proposition which will be applied in the follow in section is the case $\Omega=\Ri^d$.

\begin{cor}\label{csa2.2}
If  $\Omega=\Ri^d$ then $H=-\divv(C\nabla)$ is self-adjoint.
\end{cor}
\proof\
It follows from (\ref{esa2.3}) that
$\lambda \,\|\eta\varphi\|_2^2\leq (\varphi,\Gamma_{\!\!c}(\eta)\varphi)\leq \nu\,(\varphi, |\nabla\!\eta|^2\varphi)$
for all $\eta\in W^{1,\infty}_c(\Ri^d)$.
Now replacing $\eta$ by $\rho_n$ as above and taking the limit $n\to\infty$ leads to the conclusion that 
 $ \|\rho_n\varphi\|_2^2\to  \|\varphi\|_2^2$ and $\|\nabla\!\rho_n\|_\infty^2\to0$.
 Therefore $\varphi=0$ and $H$ is self-adjoint.
 \hfill$\Box$

\bigskip

It is more difficult to utilize this technique if $\Omega$ has a boundary $\Gamma$ but one can modify the argument slightly
to conclude that $H$ is self-adjoint if the coefficients of $H$ are sufficiently degenerate at $\Gamma$.

\begin{cor}\label{csa2.3}
If  $\Omega$  is a general domain in $\Ri^d$ but there is an $r>0$ and $\delta\geq 2$ such that $C\leq \nu\,d_\Gamma^{\,\delta}I$ on $\Gamma_{\!\!r}$ 
then $H=-\divv(C\nabla)$ is self-adjoint.
\end{cor}
\proof\
Define $\xi_n\in W^{1,\infty}(0,\infty)$ by $\xi_n(t)=0$ if $t<1/n$, $\xi_n(t)=1$ if $t>1$ and $\xi_n(t)=\log(nt)/\log n$ if $1/n\leq t\leq 1$.
Then set $\eta_n=\chi\,(\xi_n\circ(r^{-1}d_\Gamma))$ where $\chi\in C_c^\infty(\Ri^d)$ with $D=(\Omega\cap  \supp\chi)$ non-empty.
Thus the $\eta_n\in W^{1,\infty}_c(\Omega)$ have  support in $D\cap (\Omega_{r/n})$.
It follows immediately that $\lim_{n\to\infty}\|\eta_n\psi-\chi\psi\|_2=0$ for all $\psi\in L_2(\Omega)$.
Moreover, by the Cauchy--Schwarz inequality,
\[
 (\psi,\Gamma_{\!\!c}(\eta_n)\psi)
\leq (1+\varepsilon)\,(\psi,\Gamma_{\!\!c}(\chi)|(\xi_n\circ(r^{-1}d_\Gamma))|^2\psi)
+(1+\varepsilon^{-1})\,(\psi,\chi^2\,\Gamma_{\!\!c}(\xi_n\circ(r^{-1}d_\Gamma))\psi)
\]
for all $\varepsilon>0$.
But $\supp\Gamma_{\!\!c}(\xi_n\circ(r^{-1}d_\Gamma))\subseteq{\overline \Gamma}_{\!\!r}$
and $\Gamma_{\!\!c}(\xi_n\circ(r^{-1}d_\Gamma))\leq \nu \,(r^{-1}d_\Gamma)^{\delta-2}(\log n)^{-2}\leq \nu\,(\log n)^{-2}$
on its support.
Therefore
\[
\limsup_{n\to\infty}(\psi,\Gamma_{\!\!c}(\eta_n)\psi)\leq (1+\varepsilon)\,(\psi,\Gamma_{\!\!c}(\chi)|\psi)\leq \nu\,(\psi,|\nabla\!\chi|^2\psi)
\;.
\]
Now one can replace $\eta$ by $\eta_n$ in (\ref{esa2.3}) and take the limit $n\to\infty$ followed by the limit $\varepsilon\to0$ to conclude that
if $\varphi\in D(H^*)$ and $(\lambda I+H^*)\varphi=0$ then 
\[
\lambda\,\|\chi\psi\|_2^2\leq  \nu\,(\varphi,|\nabla\!\chi|^2\varphi)
\]
for all $\chi\in C_c^\infty(\Ri^d)$.
This effectively reduces the problem to the $\Ri^d$-case covered by the previous corollary.
One again deduces that $\varphi=0$ by constructing a sequence of $\chi_n$ such that $\|\chi_n\varphi\|_2^2\to\|\varphi\|_2^2$ and 
$(\varphi,|\nabla\!\chi_n|^2\varphi)\to0$ as $n\to\infty$.
\hfill$\Box$

\bigskip

There are two distinct but related problems that occur if one tries to apply the foregoing arguments to less degenerate situations, e.g.\ to operators with $C\sim d_\Gamma^{\,\delta}I$ near the boundary with $\delta<2$.
First the criterion $\ker(\lambda I+H^*)\varphi=\{0\}$ for self-adjointness is clearly a global property but self-adjointness should be determined by boundary behaviour.
Therefore one needs to reformulate the criterion appropriately.
Secondly, the inequality (\ref{esa2.3}) is not sufficiently sensitive to the boundary behaviour.
This problem arises since we totally discarded the term $h(\eta\varphi)$ in the identity (\ref{esa2.2}).
Therefore one needs to exploit more detailed properties of $h$ near the boundary.
If $\delta\geq2$ then the boundary is essentially inaccessible to the related diffusion process and this explains why the conclusion of Corollary~\ref{csa2.3} is independent of the details of the boundary.
An alternative expression of this inaccessibility is that $\Omega$ equipped with the Riemannian metric $ds^2=d_\Gamma^{\,\delta}\,dx^2$ is complete for all $\delta\geq 2$.
Thus the corresponding Riemannian distance to the boundary is infinity.

\section{Boundary estimates}\label{S3}

In this section we examine two more preparatory topics.
First we show that the self-adjointness  property $H=H^*$
 can be verified in two steps, an interior estimate and a boundary estimate.
A similar approach was taken in \cite{Rob15} for the verification of the alternative self-adjointness criterion $H=H_{\!F}$ and the following discussion
relies partly on the results in Section~2.1 of the previous paper.
Secondly, we discuss extensions of Hardy inequalities on a boundary layer to weak Hardy inequalities on the whole domain.

The first step in establishing that  $H=H^*$  is to verify the property on  the interior sets $\Omega_r=\{x\in\Omega: d_\Gamma(x)>r\}$.

\begin{lemma}\label{lsa3.1} Assume $\supp\varphi\subseteq \Omega_r$ for some $r>0$.
Then $\varphi\in D(H^*)$ if and only if $\varphi\in D(H)$.
Moreover, if these conditions are satisfied then $H^*\varphi=H\varphi$.
\end{lemma}
\proof\
First assume $\varphi\in D(H^*)$ with  $\supp\varphi\subseteq \Omega_r$.
Then fix $s,t>0$ such that $s<t<r$.
Secondly, choose a $\xi\in W^{2,\infty}(\Ri^d)$ with the properties $0\leq\xi\leq1$, $\supp\xi\subseteq \Omega_s$ and $\xi=1$ on $\Omega_t$.
Then define the (closed) operator $H_\xi=-\divv(C_\xi\nabla)$ on $L_2(\Ri^d)$ where the matrix $C_\xi$ is given by $C_\xi=\xi\,C+(1-\xi)\,\mu I$ with $\mu>0$.
Since $C>0$ on $\Omega$ and $\mu>0$ it follows that $C_\xi>0$ on $\Ri^d$.
Moreover, the coefficients of $C_\xi$ are locally Lipschitz.
Therefore the operator $H_\xi$  on $L_2(\Ri^d)$ is self-adjoint by Corollary~\ref{csa2.2}.

It follows from this construction that  $(H^*\varphi, \psi)=(\varphi, H\psi)=(\varphi, H_\xi\psi)$ for all $\psi\in C_c^\infty(\Omega_t)$.
But if $\psi\in C_c^\infty(\Omega\backslash \Omega_r)$ the relation follows by locality of $H$ and $H_\xi$ since this implies that all terms are identically zero.
Therefore the identity is valid for all $\psi\in C_c^\infty(\Omega)$ by decomposition.
Consequently $|(\varphi, H_\xi\psi)|\leq \|H^*\varphi\|_2\,\|\psi\|_2$.
Since $H_\xi$ is self-adjoint it follows that $\varphi\in D(H_\xi)$ and $H_\xi\varphi=H^*\varphi$.
But it follows from the proof of Proposition~2.2 in \cite{Rob15} that if $\supp\varphi\subseteq \Omega_r$ then $\varphi\in D(H_\xi)$ if and only if $\varphi\in D(H)$
and in this case $H_\xi\varphi=H\varphi$.
Combining these statements one concludes that $\varphi\in D(H)$ and $H\varphi=H^*\varphi$.

Conversely, if $\varphi\in D(H)$ with $\supp\varphi\subseteq \Omega_r$ then $\varphi\in D(H_\xi)$ and $H_\xi\varphi=H\varphi$ by Proposition~2.2 in \cite{Rob15}.
Then $(\varphi, H\psi)=(\varphi,H_\xi\psi)$ for all $\psi\in C_c^\infty(\Omega)$ and $|(\varphi, H\psi)|\leq \|H_\xi\varphi\|_2\,\|\psi\|_2$.
Therefore $\varphi\in D(H^*)$ and $H^*\varphi=H_\xi\varphi=H\varphi$.
\hfill$\Box$

\bigskip

Lemma~\ref{lsa3.1} now  allows one to characterize the self-adjointness of $H$ by its boundary behaviour.
This follows from the next proposition since $H$ is self-adjoint if and only if $R(\lambda I+H)=L_2(\Omega)$.

\begin{prop}\label{psa3.1}
If the range condition $(\varphi, (\lambda I+H)\psi)=0$ for a $\lambda>0$ and all $\psi\in C_c^\infty(\Omega)$ implies that $\varphi=0$ on a boundary layer 
$\Gamma_{\!\!r}$ then $H$ is self-adjoint.
\end{prop}
\proof\
The range condition implies that  $|(\varphi, H\psi)|=\lambda|(\varphi, \psi)|\leq \lambda\,\|\varphi\|_2\,\|\psi\|_2$ for all $\psi\in C_c^\infty(\Omega)$.
Therefore $\varphi\in D(H^*)$.
But by assumption this condition also implies that $\varphi=0$ on $\Gamma_{\!\!r}$.
Therefore $\supp \varphi\subseteq \Omega_s$ for some $s\in\langle0,r\rangle$.
Then $\varphi\in D(H)$ and $H\varphi=H^*\varphi$ by Lemma~\ref{lsa3.1}.
Consequently, $\varphi\in D(h)$ and 
\[
h(\varphi)=(\varphi, H\varphi)=(H^*\varphi,\varphi)=-\lambda\,\|\varphi\|_2^2
\;.
\]
Since $h(\varphi)\geq 0$ and $\lambda>0$ it follows immediately that $\|\varphi\|_2=0$ and $\varphi=0$.
\hfill$\Box$

\bigskip

Proposition~\ref{psa3.1} reduces the proof of self-adjointness of $H$ to the verification that all the elements in the kernel  of $\lambda I+H^*$ vanish on some thin
boundary layer $\Gamma_{\!\!r}$.
This will be achieved with a stronger version of the inequality (\ref{esa2.3}) used in Section~\ref{S2} which follows from the introduction of a Hardy-type lower bound on the form $h$.

The general Hardy inequality for second-order operators can be expressed as
\[
h(\varphi)\geq \|\chi\varphi\|_2^2
\]
for all $\varphi \in D(h)$ with $\chi$ a positive multiplier on $D(h)$.
Thus the weighted Hardy inequality~(\ref{esa1.2}) corresponds to the form of the operator $-\sum^d_{k=1}\partial_k\,d_\Gamma^{\,\delta}\,\partial_k$ with 
$\chi=b_{\delta,r}\,d_\Gamma^{\,\delta/2-1}\!$.
Inequalities of this type are known on  a wide variety of domains
(see \cite{BEL} and the references therein for background)  but are also known to fail in quite simple situations.
For example, the weighted Hardy inequality~(\ref{esa1.2})  fails if 
 $\Omega$ is a  unit   ball and  $\delta>1$ although it is valid 
 for all $\psi$ with support in each  boundary layer $\Gamma_{\!\!r}$ with $r<1$ (see Example~\ref{exsa5.4} below).
The boundary Hardy inequality is also valid on  thin boundary layers for domains with a $C^2$-boundary (see \cite{Rob15}, Proposition~2.9).  
In this situation the Hardy constant $b_\delta=(\delta-1)/2$.
It is, however, of greater interest that boundary estimates of this type are also valid for domains with very rough boundaries,  
e.g.\ boundaries of a fractal nature \cite{KZ} \cite{Leh3}.
We will discuss this in Section~\ref{S6}   but for the present these observations motivate the examination of this restricted form of the Hardy inequality.

The following lemma establishes that the general Hardy inequality  on $\Gamma_{\!\!r}$ extends to a weaker form of the inequality on $\Omega$.
This will be of utility in Section~\ref{S4}.

\begin{lemma}\label{lsa3.2}
Fix $r>0$ and $\mu\geq0$. Then assume there is a positive  $\chi\in\bigcap_{s>0}L_\infty(\Omega_s)$ 
 such that 
\begin{equation}
\mu\,\|\psi\|_2^2+h(\psi)\geq \|\chi\psi\|_2^2
\label{esa3.1}
\end{equation}
for all $\psi\in C_c^1(\Gamma_{\!\!r})$.

\smallskip

It follows that for each $\varepsilon>0$ there is a $\lambda_{r,\varepsilon}>0$ such that 
\begin{equation}
\lambda_{r,\varepsilon}\|\psi\|_2^2+h(\psi)\geq (1-\varepsilon)\|\chi\psi\|_2^2
\label{esa3.2}
\end{equation}
for all 
$\psi\in D(h)$.
\end{lemma}
\proof\
The $\mu$ plays no essential role in the proof of the lemma.
Its presence only changes the value of the resulting $\lambda_{r,\varepsilon}$. 
Therefore we assume in the following argument that $\mu=0$ although in Section~\ref{S5} we  use the result  with $\mu>0$.

First fix $\xi\in C^1(\Omega)$ with $0\leq\xi\leq1$, $\xi=1$ on $\Gamma_{\!\!s}$ for some $s\in\langle0,r\rangle$ and $\xi=0$ on $\Omega_r$.
Then for each $\psi\in C_c^1(\Omega)$ one has $\xi\psi\in C_c^1(\Gamma_{\!\!r})$ and 
\[
\Gamma_{\!\!c}(\xi\psi)=\xi^2\,\Gamma_{\!\!c}(\psi)+2\,\xi\psi \,\Gamma_{\!\!c}(\xi,\psi)+\psi^2\,\Gamma_{\!\!c}(\xi)
\]
by the Leibniz rule.
Hence for each $\varepsilon>0$ one has
\[
\Gamma_{\!\!c}(\xi\psi)\leq (1+\varepsilon)\,\xi^2\,\Gamma_{\!\!c}(\psi)+(1+\varepsilon^{-1})\,\psi^2\,\Gamma_{\!\!c}(\xi)
\]
by the Cauchy-Schwarz inequality.
Therefore by integration and rearrangement
\[
h(\psi)\geq \int_\Omega \xi^2\,\Gamma_{\!\!c}(\psi)\geq (1+\varepsilon)^{-1}\,h(\xi\psi)-\varepsilon^{-1}\int_{\Gamma_{\!\!r}\backslash \Gamma_{\!\!s}}\psi^2\,
\Gamma_{\!\!c}(\xi)
\]
since $\supp\Gamma_{\!\!c}(\xi)\subseteq \Gamma_{\!\!r}\backslash \Gamma_{\!\!s}$.
Hence 
\[
\varepsilon^{-1}\lambda_{r,s}\,\|\psi\|_2^2+h(\psi)\geq (1+\varepsilon)^{-1}\,h(\xi\psi)
\]
with $\lambda_{r,s}$ the $L_\infty$-norm of $\Gamma_{\!\!c}(\xi)$.
Then it follows from the assumption (\ref{esa3.1}), with $\mu=0$ and $\psi$ replaced by $\xi\psi$, that 
\[
h(\xi\psi)\geq \|\chi\xi\psi\|_2^2\geq \|\chi_s\psi\|_2^2
\]
where $\chi_s$ denotes the restriction of $\chi$ to $\Gamma_{\!\!s}$.
But 
\[
 \|\chi_s\psi\|_2^2= \|\chi\psi\|_2^2-\int_{\Gamma_{\!\!s}^c}|\chi|^2|\psi|^2\geq \|\chi\psi\|_2^2-\lambda_s\|\psi\|_2^2
 \]
 where $\lambda_s=\sup\{\chi(x):d_\Gamma(x)>s\}$.
 Combining these estimates gives
 \[
 (\lambda_s+\varepsilon^{-1}\lambda_{r,s})\,\|\psi\|_2^2+h(\psi)\geq (1+\varepsilon^{-1})\,\|\chi\psi\|_2^2
 \]
 for all $\varepsilon>0$ and all $\psi\in C_c^1(\Omega)$.
 Finally with the replacement $\varepsilon\mapsto\varepsilon^{-1}-1$ one obtains (\ref{esa3.2}) for all $\psi\in C_c^1(\Omega)$ with 
 $\lambda_{r,\varepsilon}=\inf_{s\in\langle0,r\rangle}(\lambda_s+\varepsilon^{-1}\lambda_{r,s})$.
Then (\ref{esa3.2}) follows for all $\psi\in D(h)$
 by closure.

If $\mu>0$ then  $\lambda_{r,\varepsilon}$ is replaced  by $\mu+ \lambda_{r,\varepsilon}$.
 \hfill$\Box$

\bigskip

Note that the strong Hardy inequality (\ref{esa3.1}) only involves the restriction of $\chi$ to $\Gamma_{\!\!r}$;
 the value of $\chi$ on the interior sets $\Omega_s$ with $s>r$ is arbitrary up to the boundedness hypothesis.
This freedom pertains in the setting of the next proposition which establishes a generalization of the basic  inequality (\ref{esa2.3})
 used  to discuss the self-adjointness of operators with a degeneracy of order $\delta\geq 2$ at the boundary in Corollary~\ref{csa2.3}.

\begin{prop}\label{psa3.2}
Assume that $\chi$ satisfies the Hardy inequality $(\ref{esa3.1})$ on $\Gamma_{\!\!r}$ and  fix $\lambda_{r,\varepsilon}\geq 0$, with $\varepsilon>0$,  such that 
\begin{equation}
\lambda_{r,\varepsilon}\|\psi\|_2^2+h(\psi)\geq (1-\varepsilon)\,\|\chi\psi\|_2^2
\label{esa3.3}
\end{equation}
for all $\psi\in D(h)$.
Further fix $\lambda>\lambda_{r,\varepsilon}$ and $\varphi\in D(H^*)$ such that $(\lambda I+H^*)\varphi=0$.

\smallskip

It follows that 
\begin{equation}
(\lambda-\lambda_{r,\varepsilon})\,\|\eta\varphi\|_2^2+(1-\varepsilon)\,\|\eta\chi\varphi\|_2^2\leq (\varphi, \Gamma_{\!\!c}(\eta)\varphi)
\label{esa3.4}
\end{equation}
for all $\eta\in\bigcap_{s>0} W^{1,\infty}(\Omega_s)$.
\end{prop}
\proof\
First it follows from Lemma~\ref{lsa3.2} that for each $\varepsilon>0$ one may indeed choose $\lambda_{r,\varepsilon}$ such that (\ref{esa3.3}) is satisfied.
Then the inequality (\ref{esa3.4}) follows for all $\eta\in W^{1,\infty}_c(\Omega)$ by combination of the basic identity (\ref{esa2.2}) and the weak Hardy inequality
(\ref{esa3.3}).
It then extends to the larger set of $\eta$ by repetition of the argument used to establish Corollary~\ref{csa2.1}.
\hfill$\Box$

\bigskip

In the sequel inequality (\ref{esa3.4}) plays a  similar role  in the discussion  of self-adjointness for operators with a degeneracy $\delta<2$ 
to that played by   (\ref{esa2.3}) in the proof of self-adjointness of operators with $\delta\geq 2$ in Section~\ref{S2}.

%

\section{A prototypical theorem}\label{S4}

In this section we develop the ideas of Agmon \cite{Agm1} and Nenciu and Nenciu \cite{NeN} to derive a general self-adjointness theorem which serves as  a prototype
for a more specific result in the following section.

The discussion in Section~\ref{S2} of self-adjointness of operators with coefficients $C$  satisfying  the degeneracy condition $C\leq \nu\,d_\Gamma^{\,\delta} I$, 
with $\delta\geq2$, on a boundary layer $\Gamma_{\!\!r}$ was based on the inequality  (\ref{esa2.3}).
The proof followed by choosing a sequence $\eta\in W_c^{1,\infty}(\Omega)$ such that $\eta\to\one_\Omega$ pointwise and 
$\Gamma_{\!\!c}(\eta)\to0$.
But the first  condition indicates  that $|\nabla\eta|$ could increase as $d_\Gamma^{\,-1}$ near the boundary.
In this case  it is inevitable that 
$\Gamma_{\!\!c}(\eta)\sim d_\Gamma^{\,\delta-2}$ for small $d_\Gamma$.
Hence  if  $\delta\geq2$ then $\Gamma_{\!\!c}(\eta)$ is bounded. 
This was the key feature of the arguments in Section~\ref{S2}.
If, however, $\delta<2$ then $\Gamma_{\!\!c}(\eta)$ is unbounded and the reasoning of Section~\ref{S2} is totally  inadequate.
But the improved inequality (\ref{esa3.4})  gives a potential path to circumvent this difficulty.
The factor $\Gamma_{\!\!c}(\eta)$   on the right hand side might well diverge at the boundary but  the term on the left with the factor $\eta\chi$ could also diverge.
Therefore the idea is to choose a sequence of $\eta$ such that the two divergences cancel and do not interfere with the estimation argument.
This is the strategy of Nenciu and Nenciu \cite{NeN} in the case of bounded $\Omega$ with some basic smoothness of the boundary.
But the method also extends to the more general case of unbounded domains with rough boundaries if one has sufficient information on the Hardy inequality near the boundary.
This is illustrated by the following prototypical result.

\begin{thm}\label{tsa4.1}
Assume there is an $r>0$ such that $C\leq \nu\, d_\Gamma^{\,\delta}I$ on $\Gamma_{\!\!r}$  with $\nu>0$ and  $\delta\in[0,2\rangle$.
Moreover, assume the Hardy inequality 
\begin{equation}
h(\psi)\geq \nu\,b_{\delta,r}^{\,2}\, \|d_\Gamma^{\,\delta/2-1}\psi\|_2^2
\label{esa4.00}
\end{equation}
 is satisfied with $b_{\delta,r}>0$  for all $\psi\in D(h)$ with $\supp\psi\subseteq \Gamma_{\!\!r}$.
 Let $b_\delta$ denote the supremum over small $r$ of the  possible $b_{\delta,r}$.

\smallskip

If $(2-\delta)/2<b_\delta$  then $H$ is self-adjoint.
\end{thm}
\proof\
The Hardy inequality (\ref{esa4.00}) is analogous to the weighted inequality (\ref{esa1.2})
and the Hardy constant  $b_\delta$ is defined similarly.
The inequality corresponds to (\ref{esa3.1})   with $\mu=0$
and  $\chi= (\nu^{1/2}\,b_{\delta,r})\,d_\Gamma^{\,\delta/2-1}$. 
Thus it follows from Lemma~\ref{lsa3.2} and Proposition~\ref{psa3.2} that (\ref{esa3.4}) is satisfied with this choice of $\chi$ on $\Omega$.
Note that Lemma~\ref{lsa3.2} is applicable  since $\chi$ is positive and bounded on the interior set $\Omega_r$ for  $\delta\leq 2$.

Now the principal idea is to utilize (\ref{esa3.4})
with $\varphi\in \ker(\lambda I+H^*)$ to deduce that $\varphi=0$ on $\Gamma_{\!\!r}$.
This  suffices for self-adjointness by Proposition~\ref{psa3.1}.
The proof uses  a method  introduced by Agmon with a different aim in mind (see \cite{Agm1} Theorem~1.5). 
First one expresses   $\eta$ in the form $\eta=e^\xi\,\zeta$, with suitable support restrictions on $\xi$ and $\zeta$.
Then  (\ref{esa3.4}) gives
\[
(\lambda-\lambda_{r,\varepsilon})\,\|\chi e^\xi\zeta\varphi\|_2^2
+(1-\varepsilon)\|e^\xi\zeta\varphi\|_2^2\leq (\varphi, \Gamma_{\!\!c}(e^\xi\zeta)\varphi)
\;.
\]
Secondly, one chooses $\xi$   such that $\Gamma_{\!\!c}(\xi)\leq (1-\varepsilon)\,\chi^2$ on $\Omega$ to obtain
\begin{equation}
(\lambda-\lambda_{r,\varepsilon})\,\|e^\xi\zeta\varphi\|_2^2+(e^\xi\zeta\varphi, \Gamma_{\!\!c}(\xi) e^\xi\zeta\varphi)\leq (\varphi, \Gamma_{\!\!c}(e^\xi\zeta)\varphi)
\;.
\label{esa4.0}
\end{equation}
This can be reformulated as 
\begin{equation}
(\lambda-\lambda_{r,\varepsilon})\,\|e^\xi\zeta\varphi\|_2^2\leq (e^\xi\varphi, \Gamma_{\!\!c}(\zeta)e^\xi\varphi)+2\,(e^\xi\varphi, \Gamma_{\!\!c}(\xi,\zeta)\zeta e^{\xi}\varphi)
\label{esa4.1}
\end{equation}
by evaluating the right hand side of (\ref{esa4.0}) with the Leibniz rule.
This corresponds closely to Lemma~3.4 of \cite{NeN} although the latter lemma is expressed in a different manner.
The key point is that the second term  
on the left  hand side of (\ref{esa4.0}) is  cancelled by  the leading term in the Leibniz
expansion of  the right hand side.
Finally one replaces $\zeta$ by a sequence $\zeta_n$ where $\zeta_n\to1$ in such a way that one can control the growth of the right hand side and 
subsequently deduce that $\varphi=0$ on $\Gamma_{\!\!r}$.

\smallskip

First, define $\hat\xi$ on $\langle0,\infty\rangle$ by
\[
{\hat\xi}(t)=\log(t/(1+t))^{(2-\delta)/2} +2^{-1}\log\log((1+t)/t)
\;.
\]
(This corresponds to  the function $G$ used by  Nenciu and Nenciu \cite{NeN} in the proof of their Theorem~5.3
 but with $t$ replaced by $t/(1+t)$.)
It follows immediately that  one has 
\[
e^{2{\hat\xi}(t)}=(t/(1+t))^{(2-\delta)}\log((1+t)/t)=-(2-\delta)^{-1}(t/(1+t))^{(2-\delta)}\log(t/(1+t))^{(2-\delta)}
\;.
\]
But $t/(1+t)\in \langle0,1\rangle$ for $t\in\langle0,\infty\rangle$ and $s\in\langle0,1\rangle\mapsto -s\log s$ is both positive and bounded.
Therefore $e^{2\hat \xi}$   is uniformly bounded on $\langle0,\infty\rangle$ for  $\delta\in[0,2\rangle$.
Moreover,
\[
{\hat\xi}^{\,\prime}(t)=((2-\delta)/2t)(1/(1+t))(1-((2-\delta)\log((1+t)/t))^{-1})\leq (2-\delta)/2t
\;.
\]
Now set $\xi={\hat\xi}\circ d_\Gamma$.
It follows that 
\begin{equation}
e^{2\xi}=-(2-\delta)^{-1}(d_\Gamma/(1+d_\Gamma))^{(2-\delta)}\log(d_\Gamma/(1+d_\Gamma))^{(2-\delta)}
\label{esa4.2}
\end{equation}
and 
\begin{equation}
\Gamma_{\!\!c}(\xi)\leq \nu\,((2-\delta)/2)^2\,d_\Gamma^{\,\delta-2}
\;.
\label{esa4.3}
\end{equation}
Thus  if  $((2-\delta)/2)^2\leq (1-\varepsilon)\,b_{\delta,r}^{\,2}$ then 
 the condition $\Gamma_{\!\!c}(\xi)\leq (1-\varepsilon)\,\chi^2$ used for the cancellation in passing from (\ref{esa4.0}) to (\ref{esa4.1})
is satisfied.
Note  that  $s\in\langle0,e^{-1}\rangle\mapsto -s\log s$ is increasing.
Therefore 
\begin{equation}
e^{2\xi}\leq -d_\Gamma^{\,2-\delta}\log d_\Gamma\;\;\;\;\;\mbox{and}\;\;\;\;\;e^{2\xi}\,\Gamma_{\!\!c}(\xi)\leq -\nu\,((2-\delta)/2)^2\,\log d_\Gamma
\label{esa4.4}
\end{equation}
on $\Gamma_{\!\!r}$ for all small $r$.
In particular  the term which  cancelled in the passage from (\ref{esa4.0}) to (\ref{esa4.1}) diverges logarithmically as $d_\Gamma\to0$.
\smallskip

Secondly, let $\hat \zeta\in W^{1,\infty}(0,\infty)$  be an increasing function with $0\leq \hat\zeta\leq1$, $\hat\zeta=0$ if $t<r/2$,  $\hat \zeta=1$ if $t\geq r$ and 
$|{\hat\zeta}^{\,\prime}|\leq 2/r$.
Then set $\zeta=\hat \zeta\circ d_\Gamma$ and 
$\zeta_n={\hat\zeta}\circ(2^nd_\Gamma)$.
Thus  $\zeta=\zeta_0$ and $\zeta_n=0$ if $d_\Gamma<2^{-(n+1)}r$, $\zeta_n=1$ if $d_\Gamma\geq 2^{-n}r$ and $|\nabla \zeta_n|\leq 2^{n+1}/r$.

\smallskip

Thirdly, with these choices we examine the bound (\ref{esa4.1}).
It follows by definition  that $\zeta_n=1$ on  $B_m=\{x\in\Omega:2^{-m+1}r\leq d_\Gamma\leq r\}\subset \Gamma_{\!\!r}$  if $n>m$.
But $e^{\,\xi}$ is bounded away from zero on $B_m$ by (\ref{esa4.2}).
Therefore  there is a $b_m>0$ such that the  norm on the left hand side satisfies
\begin{equation}
\|e^{\,\xi}\zeta_n\varphi\|_2^2\geq b_m\|\one_{B_m}\varphi\|_2^2
\label{esa4.41}
\end{equation}
for all $n\geq m$.
Next consider the factor $e^{2\xi}\,\Gamma_{\!\!c}(\zeta_n)$ on the right hand side of  (\ref{esa4.1}).
It clearly  has support in the set $A_n=\{x\in\Omega: 2^{-(n+1)}r\leq d_\Gamma\leq  2^{-n}r\}$ because the function $\zeta_n$ is equal to zero $0$ if $d_\Gamma<2^{-(n+1)}r$
and to  $1$ if $d_\Gamma\geq 2^{-n}r$. 
But on $A_n$ one has 
\[
\Gamma_{\!\!c}(\zeta_n)\leq 4\,\nu\,(2^{-n}r)^{(\delta-2)}
\]
by the assumed bounds on $C$ and the definition of $\zeta_n$.
Since   $e^{2\xi}\leq -(2^{-n}r)^{(2-\delta)}\log(2^{-n}r)$ on $A_n$ by (\ref{esa4.4}) it follows that
\begin{equation}
e^{2\xi}\,\Gamma_{\!\!c}(\zeta_n)\leq -4\,\nu\,\log(2^{-n}r)\leq 8\,n\,\nu
\label{esa4.5}
\end{equation}
on $A_n$ for all large $n$.
The second factor $e^{2\xi}\,\Gamma_{\!\!c}(\xi,\zeta_n)\,\zeta_n$ on the right hand side of  (\ref{esa4.1})
also has support in $A_n$ and 
\[
|e^{2\xi}\,\Gamma_{\!\!c}(\xi,\zeta_n)\,\zeta_n|\leq (e^{2\xi}\,\Gamma_{\!\!c}(\xi))^{1/2}\,(e^{2\xi}\,\Gamma_{\!\!c}(\zeta_n))^{1/2}
\]
by the Cauchy--Schwarz inequality.
Therefore it follows from (\ref{esa4.4}) and (\ref{esa4.5}) that 
\begin{equation}
|e^{2\xi}\,\Gamma_{\!\!c}(\xi,\zeta_n)\zeta_n|\leq 8\,(2-\delta)\,n\,\nu
\label{esa4.6}
\end{equation}
on $A_n$ for all large $n$.

Finally substituting the estimates (\ref{esa4.41}), (\ref{esa4.5}) and (\ref{esa4.6}) into 
(\ref{esa4.1}) one concludes that there is an $a>0$, independent of $n$, such that 
\[
(\lambda-\lambda_{r,\varepsilon})\,b_m\,\|\one_{B_m}\varphi\|_2^2\leq a\,n\,\|\one_{A_n}\varphi\|_2^2
\]
for all large $n$, and in particular $n\geq m$.
But then for all large $N_1, N_2$ one must have
\[
(\lambda-\lambda_{r,\varepsilon})\,b_m\,\Big(\sum\nolimits^{N_2}_{n=N_1}n^{-1}\Big)\,\|\one_{B_m}\varphi\|_2^2\leq a\,\Big(\sum\nolimits^{N_2}_{n=N_1}\|\one_{A_n}\varphi\|_2^2\Big)\leq a\,\|\varphi\|_2^2
\;.
\]
Since the sum on the left diverges as $N_2\to\infty$ it follows  that  $\one_{B_m}\varphi=0$. 
As this is  valid for all $m$ it follows that $\varphi=0$ on $\Gamma_{\!\!r}$.

Therefore we have now deduced that  if $((2-\delta)/2)^2\leq (1-\varepsilon)\,b_{\delta,r}^{\,2}$, which ensures that  $\Gamma_{\!\!c}(\xi)\leq (1-\varepsilon)\,\chi^2$,  
then  $H$ is self-adjoint by Proposition~\ref{psa3.1}.
But this conclusion is valid for all small $\varepsilon$ and $r$.
Thus $H$ is self-adjoint whenever $(2-\delta)/2<b_\delta$.
\hfill$\Box$

\bigskip

As pointed out in the introduction the statement of Theorem~\ref{tsa4.1} can be strengthened if the boundary $\Gamma$ of $\Omega$ separates into a countable union of 
positively separated components $\Gamma^{(j)}$.
Then the degeneracy can vary from component to component.

\begin{cor}\label{csa4.1}
Assume $\Gamma=\bigcup_{j\geq 1}\Gamma^{(j)}$ with $d(\Gamma^{(j)},\Gamma^{(k)})\geq r_0>0$ for all $j\neq k$ .
If, for each $j$ one has  $C\leq \nu\,d_{\Gamma^{(j)}}^{\;\delta_j/2-1}I$ on  $\Gamma^{(j)}_{\!\!r }$ with $r<r_0/2$ and $\delta_j\in[0,2\rangle$ and if the Hardy inequality $(\ref{esa3.1})$ is satisfied
with $\chi_j= \nu^{1/2}\,b_{\delta_j}\,d_{\Gamma^{(j)}}^{\;\delta_j/2-1}$ on  $\Gamma^{(j)}_{\!\!r}$ where $b_{\delta_j}>(2-\delta_j)/2$ then $H$ is self-adjoint.
\end{cor}
\proof\
The proof is essentially unchanged.
First one can prove that if $\varphi\in D(H^*)$ and $(\lambda I+H^*)\varphi=0$ then $\varphi=0$ on each $\Gamma^{(j)}_{\!\!r}$ 
by repetition of  the above argument  component by component.
with $\zeta$ successively replaced by $\zeta^{(j)}=\zeta\circ d_{\Gamma^{(j)}}$.
One then establishes that $\varphi=0$ on $\Gamma_{\!\!r}=\bigcup_{j\geq1}\Gamma^{(j)}_{\!\!r}$
and the proof follows as before.
\hfill$\Box$

\section{A direct theorem}\label{S5}

Theorem~\ref{tsa4.1} established a self-adjointness criterion  from  a weak but explicit upper bound  on the coefficients $C$ and an implicit lower bound, a Hardy inequality on a  boundary layer.
In this section we show that self-adjointness also follows from the  asymptotic degeneracy condition on the coefficients $C$ at the boundary used in the earlier paper \cite{Rob15}
together with  the  weighted Hardy inequality (\ref{esa1.2}).
This form of the result is  convenient in verifying self-adjointness for particular types of domain.

Throughout the sequel the coefficients are assumed to satisfy the boundary condtion
 \begin{equation}
\textstyle{\inf_{r\in\langle0,r_0]}}\;\textstyle{\sup_{x\in\Gamma_{\!\!r}}}\|(C\,d_\Gamma^{\,-\delta})(x)-c(x)I\|=0
\label{esa5.1}
 \end{equation}
 for some $r_0>0$
 where $c$ is  a  bounded Lipschitz function satisfying  $\inf_{x\in\Gamma_{\!\!r}}c(x)\geq \mu>0$ and $\delta\geq0$.
Condition~(\ref{esa5.1})  can be interpreted in an obvious way as
\[
\textstyle{\limsup_{d_\Gamma\to0}} \;C(c\,d_\Gamma^{\,\delta}I)^{-1}=I
\;.
\]
Thus in the language of asymptotic analysis  $C$ converges to  $c\,d_\Gamma^{\,\delta} I$ as $d_\Gamma\to0$ (see \cite{DeB}).
The  parameter $\delta$ determines the order of degeneracy at the boundary and  $c$ describes  the boundary  profile of $C$. 

The comparability of $C$ and $c\,d_\Gamma^{\,\delta}I$ can be made more precise by noting
that  for each $r\in\langle0,r_0]$  there are $\sigma_{\!r}, \tau_{\!r}>0$ such that 
\begin{equation}
\sigma_{\!r}(c \,d_\Gamma^{\,\delta})(x)I\leq C(x)\leq \tau_{\!r} (c \,d_\Gamma^{\,\delta})(x)I
\label{esa5.2}
\end{equation}
for all $x\in  \Gamma_{\!\!r}$.
The earlier discussions of Markov uniqueness  in \cite{RSi4} and \cite{LR} were based on these latter conditions.
They  play a key role in the following together with the observation that the limit condition (\ref{esa5.1}) implies that 
$\sigma_{\!r},\tau_{\!r}\to1$ as $r\to0$.
In fact one may assume that $\sigma_{\!r}$ converges monotonically upward  and $\tau_{\!r}$ converges monotonically downward.

\smallskip

Secondly, we assume that the weighted Hardy inequality (\ref{esa1.2}) is valid on $\Gamma_{\!\!r}$. 
Thus for each $r\in\langle0,r_0]$ and $\delta\geq 0$ there is a $b_{\delta,r}>0$ such that
\begin{equation}
\int_{\Gamma_{\!\!r}}d_\Gamma^{\,\delta}\,|\nabla\psi|^2\geq b_{\delta, r}^{\,2}\int_{\Gamma_{\!\!r}}d_\Gamma^{\,\delta-2}\,|\psi|^2
\label{esa5.3}
\end{equation}
for all $\psi\in C_c^1(\Gamma_{\!\!r})$.
Then the Hardy constant  $b_\delta=\sup_{r>0} b_{\delta,r}$ where again  the supremum is over  all possible $b_{\delta,r}$ for which (\ref{esa5.3}) holds.
Although the  weighted Hardy inequality (\ref{esa5.3}) near the boundary is independent of $c$ it
 does lead to a weak form of the Hardy inequality for the operator $H$ or, more precisely,  for the form $h$, on the whole domain  $\Omega$.

\begin{lemma}\label{lsa5.1}
Assume the boundary condition $(\ref{esa5.1})$ and the  weighted Hardy inequality $(\ref{esa5.3})$ are valid on the boundary layer $\Gamma_{\!\!r}$. 
Then  for each $\varepsilon>0$ there is a $\lambda_{r,\varepsilon}>0$ such that
\[
\lambda_{r,\varepsilon}\,\|\psi\|_2^2+h(\psi)\geq \sigma_{\!r}(1-\varepsilon)^2\,b_{\delta, r}^{\,2}\,\|c^{1/2}d_\Gamma^{\,\delta/2-1}\psi\|_2^2
\]
for all $\psi\in C_c^1(\Omega)$.
\end{lemma}
\proof\
It follows from (\ref{esa5.2}) that 
\begin{eqnarray*}
h(\psi)\geq \sigma_{\!r}\int_{ \Gamma_{\!\!r}}c\,d_\Gamma^{\,\delta}\,|\nabla\psi|^2
=\sigma_{\!r}\int_{ \Gamma_{\!\!r}}\,d_\Gamma^{\,\delta}\,|\nabla( c^{1/2}\psi)- (\nabla c^{1/2})\psi|^2
\end{eqnarray*}
 for all $\psi\in C_c^\infty(\Gamma_{\!\!r})$.
Then  by the Cauchy--Schwarz inequality  one has for each $\varepsilon>0$ 
\[
h(\psi)\geq \sigma_{\!r}(1-\varepsilon)\int_{ \Gamma_{\!\!r}}\,d_\Gamma^{\,\delta}\,|\nabla( c^{1/2}\psi)|^2
- \sigma_{\!r}(\varepsilon^{-1}-1)\int_{ \Gamma_{\!\!r}}\,d_\Gamma^{\,\delta}\,|(\nabla c^{1/2})\psi|^2
\]
for all $\psi\in C_c^\infty(\Gamma_{\!\!r})$.
Therefore, by the weighted  Hardy inequality  (\ref{esa5.3}), one deduces that 
\[
\mu_{r,\varepsilon}\|\psi\|_2^2+h(\psi)\geq \sigma_{\!r}(1-\varepsilon)\,b_{\delta, r}^{\,2}\,\int_{ \Gamma_{\!\!r}}\,c\,d_\Gamma^{\,\delta-2}\,|\psi|^2
\]
for all $\psi\in C_c^1(\Gamma_{\!\!r})$ with $\mu_{r,\varepsilon}=\sigma_{\!r}(\varepsilon^{-1}-1)r^\delta(\|\nabla c\|_\infty^2/\lambda)$.
Finally it follows from Lemma~\ref{lsa3.2}  with $\mu=\mu_{r,\varepsilon}$ and $\chi^2= \sigma_{\!r}(1-\varepsilon)\,b_{\delta, r}^{\,2}\,c\,d_\Gamma^{\,\delta-2}$
that there is a $\lambda_{r,\varepsilon}>0$ such that 
\[
\lambda_{r,\varepsilon}\,\|\psi\|_2^2+h(\psi)\geq \sigma_{\!r}(1-\varepsilon)^2\,b_{\delta, r}^{\,2}\,\|c^{1/2}d_\Gamma^{\,\delta/2-1}\psi\|_2^2
\]
for all $\psi\in C_c^1(\Omega)$.
\hfill$\Box$

\bigskip

Now one has the direct version of Theorem~\ref{tsa4.1}.

\begin{thm}\label{tsa5.1}
Assume the coefficient matrix $C$ satisfies  the boundary condition $(\ref{esa5.1})$ 
and the  weighted boundary Hardy inequality $(\ref{esa5.3})$  inequality is valid  with $\delta\in[0,2\rangle$.

\smallskip

If $(2-\delta)/2<b_\delta$ then $H$ is self-adjoint.
\end{thm}
\proof\
The proof is very similar to the proof of Theorem~\ref{tsa4.1} but with some small changes.
First the upper bound on the coefficient matrix $C$ is now replaced by the upper bound in (\ref{esa5.2}).
Secondly,  one observes that since the  weighted Hardy inequality (\ref{esa5.3}) is valid the weak Hardy inequality of Lemma~\ref{lsa5.1} is also valid.
But then one obtains an inequality  identical in form to (\ref{esa3.4}) with a slightly different choice of $\chi$.
Nevertheless the  choices of $\eta$, $\xi$ and $\zeta$ are as before. 
Now, however, to verify the cancellation condition $\Gamma_{\!\!c}(\xi)\leq(1-\varepsilon)\,\chi^2$ we have to take into account the modified bound on $C$ and the different choice of 
 $\chi$.

First,  since  the bound $C\leq \nu\,d_\Gamma^{\,\delta }I$ used previously to estimate $\Gamma_{\!\!c}(\xi)$ is now replaced by right hand bound 
$C\leq \tau_{\!r} \,c \,d_\Gamma^{\,\delta} I$ of (\ref{esa5.2}) one effectively
replaces  $\nu$  by $ \tau_{\!r} \,c$ in the upper bound on  $\Gamma_{\!\!c}(\xi)$.
Explicitly, the earlier bound $\Gamma_{\!\!c}(\xi)\leq \nu\,d_\Gamma^{\,\delta}\,|\nabla\xi|^2$ is replaced by $\Gamma_{\!\!c}(\xi)\leq \tau_{\!r}\,(c\,d_\Gamma^{\,\delta})\,|\nabla\xi|^2$
with a  similar change in the bound on $\Gamma_{\!\!c}(\zeta)$.

Secondly, the earlier argument used the Hardy inequality $h(\psi)\geq \|\chi\psi\|_2^2$  for $\psi\in C_c^\infty(\Gamma_{\!\!r})$ with the identification $\chi^2=\nu\,b_{\delta,r}^{\,2}\,d_\Gamma^{\,\delta-2}$.
This then led to the weak Hardy inequality 
\[
\lambda_{r,\varepsilon}\|\psi\|_2^2+h(\psi)\geq (1-\varepsilon)\,\|\chi\psi\|_2^2
\]
for all $\psi\in D(h)$ by Lemma~\ref{lsa3.2}.
Now, however, the weighted Hardy inequality (\ref{esa5.3}) gives the analogous weak Hardy inequality of Lemma~\ref{lsa5.1} but with
$\chi^2= \sigma_{\!r}(1-\varepsilon)\,b_{\delta, r}^{\,2}\,(c\,d_\Gamma^{\,\delta-2})$
Therefore  $\nu\,b_{\delta,r}^{\,2}\,d_\Gamma^{\,\delta-2}$ is  replaced by $\sigma_{\!r}\,(1-\varepsilon)\, b_{\delta,r}^{\,2}\,(c\,d_\Gamma^{\,\delta-2})$ in the identification of $\chi^2$.

Thirdly, 
 after these replacements 
the condition  $\Gamma_{\!\!c}(\xi)\leq(1-\varepsilon)\,\chi^2$  which previously translated to  
$\nu\,((2-\delta)/2)^2\,d_\Gamma^{\,\delta-2}\leq(1-\varepsilon)\, \nu\,b_{\delta,r}^{\,2}\,d_\Gamma^{\,\delta-2}$
is replaced by  the similar inequality $ \tau_{\!r} \, ((2-\delta)/2)^2\,(c\,d_\Gamma^{\,\delta-2})\leq (1-\varepsilon)^2\,\sigma_{\!r}\, b_{\delta,r}^{\,2}\,(c\,d_\Gamma^{\,\delta-2})$.
Thus after cancellation of the strictly positive function $c\,d_\Gamma^{\,\delta-2}$ one obtains the condition
$((2-\delta)/2)^2\leq  (1-\varepsilon)^2\,(\sigma_{\!r}/\tau_{\!r})\, b_{\delta,r}^{\,2}$ for each $r\in \langle0,r_0\rangle$ and $\varepsilon\in\langle0,1\rangle$.

Fourthly, all the arguments of the previous proof carry through with these modifications. 
Therefore one concludes that for fixed $r$ the condition 
$((2-\delta)/2)^2\leq  (1-\varepsilon)^2\,(\sigma_{\!r}/\tau_{\!r})\, b_{\delta,r}^{\,2}$
suffices for self-adjointness of $H$.

Finally one may take the essential supremum of the right hand side of this condition over $r$ followed by the limit $\varepsilon\to0$  to deduce that 
 $((2-\delta)/2)^2< b_{\delta}^{\,2}$ is sufficient for  self-adjointness.
\hfill$\Box$

\bigskip

Theorem~\ref{tsa5.1} reduces the proof of self-adjointness to verification of  the boundary Hardy inequality (\ref{esa5.3}) and calculation of  the corresponding
boundary Hardy constant $b_\delta$.
Both these problems have been resolved with additional smoothness or convexity assumptions for the boundary. 
A general result of this nature is the following.

\begin{lemma} \label{lxsa5.1} Assume there are $\beta, \gamma>0$ such that
\begin{equation}
|d_\Gamma(\nabla^2d_\Gamma)-(\beta-1)|\leq \gamma\,d_\Gamma
\label{exsa5.1}
\end{equation}
in the weak sense on a boundary layer $\Gamma_{\!\!r}$.
Then if $\delta+\beta-2>\gamma r$ the weighted Hardy inequality $(\ref{esa5.3})$ is valid on $\Gamma_{\!\!r}$
with $b_{\delta,r}=(\delta+\beta-2-\gamma r)/2$.
In particular $b_\delta=(\delta+\beta-2)/2$ under the restriction $\delta>2-\beta$.
\end{lemma}
\proof\
The proof is based on a standard argument which will be applied to other settings below.
It was used in  the proof of Proposition~2.9 in \cite{Rob15}.
Set $\chi=d_\Gamma^{\,\delta-1}(\nabla d_\Gamma)$.
Then 
\[
\divv \chi=(\delta-1+d_\Gamma (\nabla^2 d_\Gamma))\,d_\Gamma^{\,\delta-2}\geq (\delta+\beta-2-\gamma\,r)\,d_\Gamma^{\,\delta-2}
\]
on $\Gamma_{\!\!r}$ where we have used (\ref{exsa5.1}).
Therefore if $\delta+\beta-2-\gamma\,r>0$ one has 
\begin{eqnarray*}
(\delta+\beta-2-\gamma\,r)\int_{\Gamma_{\!\!r}}d_\Gamma^{\,\delta-2}\psi^2&\leq &\int_{\Gamma_{\!\!r}}(\divv\chi)\psi^2\\[5pt]
&\leq&2\int_{\Gamma_{\!\!r}}|\chi.\nabla\psi|\,|\psi|
\leq 2\Big(\int_{\Gamma_{\!\!r}}d_\Gamma^{\,\delta-2}\psi^2\Big)^{1/2}\,\Big(\int_{\Gamma_{\!\!r}}d_\Gamma^{\,\delta}(\nabla\psi)^2\Big)^{1/2}
\end{eqnarray*}
for all $\psi\in C_c^\infty(\Gamma_{\!\!r})$.
Then the Hardy inequality follows by squaring the last inequality and dividing out the common factor.
The identification of $b_\delta$ is immediate.
\hfill$\Box$

\bigskip

Condition (\ref{exsa5.1}) was derived for $C^2$-domains and subdomains by Brusentsev \cite{Brus}, Section~6,  and also exploited
by Filippas, Maz'ya and Tertikas \cite{FMT1}, Section~4,  in their analysis of Hardy--Sobolev inequalities.
It was  also used to establish Theorem~5.3 in \cite{NeN1} and Theorem~3.1 in \cite{Rob15}.
Combination of Theorem~\ref{tsa5.1} and Lemma~\ref{lxsa5.1}, however,    yield a stronger version 
of this latter result.
\begin{cor}\label{cexsa5.2}
Assume that the coefficients $C$ of  $H$ satisfy the boundary  condition $(\ref{esa5.1})$.
Further assume  either $\Omega$ is  a $C^2$-domain in $\Ri^d$, or $\Omega=\Ri^d\backslash\{0\}$, or $\Omega=\Ri^d\backslash\overline \Pi$
with $\Pi$ a  $C^2$-domain in the subspace $\Ri^s$.

\smallskip
It follows that $H$ is self-adjoint whenever $\delta>2-(d-d_{\!H})/2$.
\end{cor}
\proof\ The conclusion follows since in each case $d_\Gamma$ satisfies (\ref{exsa5.1}) with $\beta=d-d_{\!H}$ (see \cite{Rob15} Subsection~2.3).
Therefore $b_\delta=(d-d_{\!H}+\delta-2)/2$, by Lemma~\ref{lxsa5.1}.
Then the sufficient  condition $b_\delta>(2-\delta)/2$ for self-adjointness of Theorem~\ref{tsa5.1} 
is equivalent to $\delta>2-(d-d_{\!H})/2$.
\hfill$\Box$

\bigskip
The earlier version of this result,  Theorem~3.1 in \cite{Rob15},  required some explicit bounds on the derivatives of the coefficients
of the operator on the boundary layer $\Gamma_{\!\!r}$.
In the above corollary  no such constraint is necessary.

\smallskip

One can also apply Theorem~\ref{tsa5.1} to domains which are the complement of convex sets.

\begin{prop}\label{pexsa5.3}
Assume that the coefficients $C$ of  $H$ satisfy the boundary  condition $(\ref{esa5.1})$.
Further assume that $\Omega=\Ri^d\backslash K$  where $K$ is a non-empty closed convex set.

\smallskip
It follows that $H$ is self-adjoint whenever $\delta>2-(d-d_{\!H})/2$.
\end{prop}
\proof\
First consider the case $\dim(K)=d$. Then $d_{\!H}=d-1$ and the condition to establish is $\delta>3/2$.
But  it follows from the convexity of $K$ that $d_\Gamma$ is convex on open convex subsets of $\Omega$ and $\nabla^2d_\Gamma\geq0$
(see, for example, \cite{Rob13} Proposition~2.1).
Therefore repeating the argument in the proof of Lemma~\ref{lxsa5.1} gives $\divv\chi\geq (\delta-1)/2$.
Then the weighted Hardy inequality is valid for $\delta>1$ with $b_{\delta,r}=b_\delta=(\delta-1)/2$
and  the sufficiency condition $b_\delta>(2-\delta)/2$ is equivalent to $\delta>3/2$.

Secondly, assume $\dim(K)=s$ with $s\in\{1,\ldots, d-1\}$.
Then $d_{\!H}=s$.
Now one can factorize $\Ri^d=\Ri^{s}\times \Ri^{d-s}$ such that $K$ is a closed convex subset of $\Ri^s$.
Choose coordinates $x=(y,z)$  with $y\in \Ri^s$ and $z\in \Ri^{d-s}$. Then for $x\in \Omega$ one has
$d_\Gamma(x)=(d_K(y)^2+d_A(z)^2)^{1/2}$ where $d_K$ is the Euclidean distance to $K$ in $\Ri^s\backslash K$
and $d_A(z)=|z|$.
Then it follows as in \cite{Rob15} that 
\[
d_\Gamma(\nabla_{\!\!x}^2d_\Gamma)=d_K(\nabla_{\!\!y}^2d_K)+(d-s-1)
\;.
\]
But now $\nabla_{\!y}^2d_K\geq0$ and the argument of Lemma~\ref{lxsa5.1} gives 
$\divv\chi\geq (\delta+d-s-2)/2$.
So one now has the  weighted Hardy inequality for $\delta>2-(d-d_{\!H})$ with $b_{\delta,r}=b_\delta=(d-d_{\!H}+\delta-2)/2$.
Therefore  the sufficiency condition for self-adjointness  $b_\delta>(2-\delta)/2$ of Theorem~\ref{tsa5.1} is again equivalent to $\delta>2-(d-d_{\!H})/2$.
\hfill$\Box$

\bigskip

It is to be expected that there is a similar result to Corollary~\ref{cexsa5.2} or Proposition~\ref{pexsa5.3} for convex domains.
But the relevant Hardy inequality is not known at present.
It is known that the weighted Hardy inequality (\ref{esa5.3}) is valid for convex domains, bounded or unbounded, if $\delta\in[0,1\rangle$
with $b_\delta=b_{\delta,r}=(1-\delta)/2$
(see \cite{Avk1} for the general $L_p$-case).
It is also known that it is valid for $\delta>1$ for unbounded convex domains or on thin boundary layers
for bounded convex domains but it is not known whether  the Hardy constant $b_\delta=(\delta-1)/2$.
We will return to this discussion in a broader context in the next section.
The convex situation is illustrated by the following example which also anticipates the more general results of 
Section~\ref{S6}.

\begin{exam}\label{exsa5.4}
Let $\Omega=B(0\,;1)$, the unit ball centred at the origin.
Then $d_\Gamma(x)=1-|x|$, $(\nabla d_\Gamma)(x)=-x/|x|$ and $(\nabla^2d_\Gamma)(x)=-(d-1)(1-|x|)/|x|$.
Thus if $r\in\langle0,r_0]$ with $r_0<1$ then
$0\geq d_\Gamma(\nabla^2d_\Gamma)\geq -(d-1)r/(1-r_0)$.
In particular (\ref{exsa5.1})  is satisfied with $\beta=1$ and $\gamma=(d-1)/(1-r_0)$.
Repeating the argument used to prove Corollary~\ref{cexsa5.2} one then obtains the weighted Hardy inequality (\ref{esa5.3})  on the boundary layer $\Gamma_{\!\!r}$
with $b_{\delta,r}=(\delta-1-\gamma\,r)/2$ for all $\delta>1+\gamma \,r$ and  $r\in \langle0,r_0]$.
Therefore $b_\delta\geq (\delta-1)/2$ for $\delta>1$.  
Conversely if $\delta>1$ one can construct a sequence $\psi_n\in C^1(\Gamma_{\!\!r})$ such that 
 $\int d_\Gamma^{\,\delta}\, |\nabla\psi_n|^2/\int  d_\Gamma^{\,\delta-2}\,|\psi_n|^2\to (\delta-1)^2/4$  as $n\to \infty$ (see Proposition~\ref{psa6.3}).
 Therefore $b_\delta=(\delta-1)/2$.
 
Note that the foregoing argument  requires $r<1$. 
In fact the weighted  Hardy inequality on $\Omega$ fails for $\delta>1$.
Specifically $\int_\Omega d_\Gamma^{\,\delta}\, |\nabla\psi|^2/\int_\Omega d_\Gamma^{\,\delta-2}\,|\psi|^2=\alpha^2$ if
 $\psi=d_\Gamma^{\,\alpha}$ with $\alpha\in\langle0,(\delta-1)/2\rangle$
where the upper bound on $\alpha $ ensures that $d_\Gamma^{\,\alpha+\delta/2-1}\in L_2(\Gamma_{\!\!r})$.
 Therefore
 \[
 \textstyle{\inf_{\psi\in C_c^1(\Omega)}}\int_\Omega d_\Gamma^{\,\delta}\, |\nabla\psi|^2/\int_\Omega d_\Gamma^{\,\delta-2}\,|\psi|^2=0
 \]
 and the Hardy inequality fails on $\Omega$.
\end{exam}

Self-adjointness of $H$ should follow in Corollary~\ref{cexsa5.2} and Proposition~\ref{pexsa5.3}  from the slightly more general condition  $\delta\geq 2-(d-d_{\!H})/2$ but the critical case $\delta= 2-(d-d_{\!H})/2$ does not follow from the
 foregoing arguments or from the arguments of \cite{Rob15}.
The conjecture is partially supported by the results for the $C^2$-case.
If $\Omega$ is a $C^2$-domain then $d_{\!H}=d-1$ and the sufficient condition of Corollary~\ref{cexsa5.2} is $\delta>3/2$.
It follows, however,  from Theorem~3.2 of \cite{Rob15} that  the matching condition $\delta\leq 3/2$ is necessary for self-adjointness.
Note that the proof of this latter result does use a bound on the derivatives of the coefficients through a bound on $|\divv (C d_\Gamma^{\,-\delta})|$ and it is not clear whether
it still holds in the more general setting of the current paper.

Finally we note that Theorem~\ref{tsa5.1}, Corollary~\ref{cexsa5.2} and Proposition~\ref{pexsa5.3} can again  be used as building blocks for the consideration of more complicated domains whose boundaries decompose into separated 
components.
In particular these results extend 
 to the case that $\Gamma$ consists of a sum of positively separated components  in exactly the same manner that Theorem~\ref{tsa4.1} extended to  Corollary~\ref{csa4.1}. 
The value of the parameter $\delta$ can vary component by component in the degeneracy condition (\ref{esa5.1}) and the Hardy inequality (\ref{esa5.3}).
For example, if $\Omega=\Ri^d\backslash S$ with $S$ a countable set of positively separated points then the condition $\delta>2-d/2$ is sufficient for self-adjointness.

\section{Rough boundaries}\label{S6}

Theorems~\ref{tsa4.1} and \ref{tsa5.1} demonstrate that  the weighted Hardy inequality at the boundary is  the key property underlying self-adjointness.
Moreover, Corollary~\ref{cexsa5.2} and Proposition~\ref{pexsa5.3} give a variety of cases in which the weighted Hardy inequality is valid and the Hardy constant can be
explicitly calculated.
This has the consequence of identifying the critical degeneracy for self-adjointness as $\delta_c=2-(d-d_{\!H})/2$.
But all these cases rely on some smoothness or convexity properties of the boundary of the domain.
In this section we examine the situation of rough boundaries which is much more opaque.

There are several  possible interpretations of  domains with rough boundary.
Two established  notions  are John domains and  the subclass of uniform domains.
For simplicity we will consider uniform domains although many of the following conclusions follow for John domains with some local uniformity.
In addition there is the concept of Ahlfors regularity which, despite its name, describes boundary sets of 
a very irregular nature, e.g.\ self-similar fractals.
The Ahlfors property was used in \cite{LR} to characterize Markov uniqueness.
We begin by summarizing 
some relevant definitions and results.

First,   $\Omega$ is defined to be a uniform  domain if there is a $\sigma\geq1$ and  for each pair of points $x,y\in \Omega$  a rectifiable  curve 
$\gamma\colon[0,l]\to\Omega$, parametrized by arc length,  such that $\gamma(0)=x$, $\gamma(l)=y$  with arc length
$l(\gamma(x\,;y))\leq \sigma\,|x-y|$ and  $d_\Gamma(\gamma(t))\geq \sigma^{-1}\,(t\wedge(l-t))$ for all $t\in [0,l]$.
Uniform domains were introduced by Martio and Sarvas \cite{MarS} as a special subclass of domains studied earlier by John \cite{John} in which the bound 
on the length of the curve $\gamma$ is omitted.
In fact these authors only examined bounded domains and  the extension to
unbounded domains was given subsequently by V{\"a}is{\"a}l{\"a} \cite{Vai1} (see also \cite{Leh3}, Section~4).
In the  case of bounded John domains the definition can be simplified.
Then it suffices that there is a preferred point $x_c\in\Omega$, the centre point, which can be connected to every other $x\in\Omega$ by a curve with the foregoing properties.
In the sequel we are interested in boundary properties of the  domains. 
Then  the uniformity property  is assumed to be valid  for all pairs of points $x,y$ in the boundary layer $\Gamma_{\!\!r}$ but    the curve $\gamma$ joining the points,
which   lies in $\Omega$,  is  not  constrained to $\Gamma_{\!\!r}$.

Secondly,  let  $B(x_0\,;r_0)$ denote  the open  Euclidean ball centred at $x_0$ with radius $r_0$.
Then the boundary $\Gamma$ is defined to be  Ahlfors $s$-regular
if  there is a regular Borel measure $\mu$ on $\Gamma$ and an $ s>0$ such that 
for each subset $A=\Gamma \cap B(x_0\,;r_0)$, with $x_0\in\Gamma$ and $r_0>0$, 
 one can choose
$a>0$ so that 
\begin{equation} 
a^{-1}\,r^s \leq \mu(A\cap B(x\,;r)) \leq a\, r^s
\label{euni1.3}
\end{equation}
for all $x\in A$ and $r\in \langle0, 2r_0\rangle$.
This is a locally uniform version of the Ahlfors regularity  property used in the theory of metric spaces 
(see, for example, the monographs  \cite{DaS, Semm, Hei, MaT}). 
It implies that $\mu$ and the Hausdorff measure $\ch^s$ on $\Gamma$ are locally equivalent and $s=d_{\!H}$, 
the Hausdorff dimension of~$\Gamma$.
The Ahlfors property  implies  that $\Gamma$ is regular in the sense that each of the subsets 
$\Gamma_{\!\!x,r}=\Gamma \cap B(x\,;r)$ with $x\in\Gamma$ has  Hausdorff dimension $s$
but $\Gamma$ could have a wildly irregular fractal nature.

As a prelude to the discussion of self-adjointness we begin with a characterization of Markov uniqueness of operators  satisfying the boundary condition (\ref{esa5.1}).
Recall that $H$ is defined to be Markov unique if it has a unique self-adjoint extension on $L_2(\Omega)$ which generates
a positive semigroup, i.e.\ a semigroup which maps positive functions into positive functions.

\begin{prop}\label{psa6.1}
Assume that $\Omega$ is a uniform domain whose boundary $\Gamma$ is Ahlfors $s$-regular.
Further assume the coefficients of $H$ satisfy the boundary condition $(\ref{esa5.1})$.

Then the following conditions are equivalent:
\begin{tabel}
\item\label{psa6.1-1}
$H$ is Markov unique,
\item\label{psa6.1-2}
$\delta\geq 2-(d-s)\;\;(\,= 2-(d-d_{\!H}))$.
\end{tabel}
\end{prop}

The proof of the proposition is a repeat of the proof of Theorem~1.1 in \cite{LR}.
In the latter reference it was assumed that $C$ satisfied bounds
$a\,d_{\Gamma}^{\,\delta}I\leq C\leq b\,d_{\Gamma}^{\,\delta}I$ on a boundary layer $\Gamma_{\!\!r}$
with $a, b>0$  constant.
But these bounds  follow from (\ref{esa5.2}) since $c$ is bounded, positive and bounded away from zero on $\Gamma_{\!\!r}$.
Then the proof that Condition~\ref{psa6.1-2} implies Condition~\ref{psa6.1-1}  is by a capacity estimate which 
relies solely on the Ahlfors regularity. 
It does not require
the uniform domain property.
The proof of the converse does require the uniform property if $d_{\!H}\in [d-1,d\rangle$ but only in one small bounded neighbourhood of $\Gamma$.
We refer to \cite{LR} for details.

\medskip

The relevance of the proposition for self-adjointness is the following.

\begin{cor}\label{csa6.1}
Under the assumptions of Proposition~$\ref{psa6.1}$ the condition $\delta\geq 2-(d-d_{\!H})$ is necessary for self-adjointness.
\end{cor}
\proof\ 
If $\delta<2-(d-d_{\!H})$ then $H$ is not Markov unique. Hence it is not self-adjoint. 
\hfill$\Box$

\bigskip

Now we have the following crucial existence result for the weighted Hardy inequality on John domains.

\begin{prop}\label{psa6.2} $(${\rm Lehrb{\"a}ck}$)$
Assume that $\Omega$ is a John domain whose boundary $\Gamma$ is Ahlfors $s$-regular.
Then for each $\delta>2-(d-d_{\!H})$
  there are $ r_0>0$ and for each $r\in\langle0,r_0\rangle$ a $b_{\delta,r}>0$ such that
\begin{equation}
\|d_\Gamma^{\,\delta/2}\,\nabla\psi\|_2\geq b_{\delta,r}\,\|d_\Gamma^{\,\delta/2-1}\psi\|_2
\label{esa6.1}
\end{equation}
for all $\psi\in C_c^1(\Gamma_{\!\!r})$.
\end{prop}
\proof\
First note that (\ref{esa6.1}) is just an operator version of the weighted Hardy inequality (\ref{esa1.2}).
Thus if  $\Omega$ is unbounded it  follows for all $C_c^1(\Omega)$ by Theorem~1.3 of \cite{Leh3}
but  with $s$ the Aikawa dimension of the boundary which is larger than the  Hausdorff dimension in general.
It follows, however,  from the Ahlfors regularity of the boundary that the Hausdorff dimension and the Aikawa dimension are equal
(see Lemma~2.1 of \cite{Leh3}).
Therefore the statement of the proposition follows by restriction to $\psi$ with support in $\Gamma_{\!\!r}$.

If $\Omega$ is bounded then the weighted Hardy inequality is not generally valid on  the whole domain.
For example,  it fails on the unit ball (see Example~\ref{exsa5.4}).
Nevertheless if $r_0$ is sufficiently small the arguments of \cite{Leh3} establish the Hardy inequality on the boundary layer.
In fact the argument simplifies considerably \cite{Leh5}.
The idea is to prove that if  $r\in\langle0,r_0\rangle$ then there is a  $B_{\delta,r}>0$ such that
\[
B_{\delta,r}^{\,2}\int_{\Gamma_{\!\!r}}d_\Gamma^{\,\delta-2}\,|\psi|^2
\leq \int_{\Gamma_{\!\!r}}d_\Gamma^{\,\delta-1}\,|\psi|\,|\nabla\psi|+\int_{\Gamma_{\!\!r}}d_\Gamma^{\,\delta}\,|\nabla\psi|^2
\]
for all $\psi\in C_c^1(\Gamma_{\!\!r})$.
Then the weighted Hardy inequality follows by applying the Cauchy-Schwarz inequality to the first term on the right hand side and rearranging.
But  the latter inequality is established by an argument involving a covering of the John domain $\Omega$ by Whitney cubes satisfying the Boman chain condition
(see, for example, \cite{BKL} \cite{HK}).
Since $\Omega$ is bounded it suffices to consider chains which begin with a cube $Q_c$ which contains the centre point $x_c$.
Then, however, one can suppose that $r$ is sufficiently small that $Q_c$ has an empty intersection with $\Gamma_{\!\!r}$.
Therefore the average of each $\psi\in  C_c^1(\Gamma_{\!\!r})$ over $Q_c$ is equal to zero 
and  the main assumption (8) in Theorem~3.1 of \cite{Leh3} is automatically satisfied.
Then repeating the proof of the Theorem~3.1, with $\Omega=\Omega'=\Gamma_{\!\!r}$, one obtains the desired inequality.
\hfill$\Box$

\bigskip

Although the principal interest in Proposition~\ref{psa6.2} is its applicability to domains with rough or uniformly disconnected boundaries
it also gives information for simpler cases such as convex domains.
In particular it establishes that the weighted Hardy inequality (\ref{esa6.1}) is satisfied on a boundary layer for a general convex domain
and for all $\delta>1$.
This is consistent with the behaviour exhibited for the unit ball in Example~\ref{exsa5.4}.
In this particular example one has $b_\delta=(\delta-1)/2$ and it could well be that this remains the case for  the Hardy constant and  a general convex domain.
This is not currently known.
It would follow from an estimate $0\geq d_\Gamma(\nabla^2d_\Gamma)\geq \gamma r$ on $\Gamma_{\!\!r}$.

\smallskip

The Hardy inequality (\ref{esa6.1}) of Proposition~\ref{psa6.2} is the main ingredient in Theorem~\ref{tsa5.1}.
Therefore if  $\delta\in \langle2-(d-d_{\!H}), 2\rangle$ and  the coefficients of $H$ satisfy the boundary condition (\ref{esa5.1}) it follows that 
the condition $b_\delta>(2-\delta)/2$ is sufficient for self-adjointness of $H$.
But for this criterion to be of utility it is necessary to have further information on the Hardy constant $b_\delta$.
It follows from  \cite{Rob15} that  $b_\delta=(d-d_{\!H}+\delta-2)/2$  for  the domains covered by Corollary~\ref{cexsa5.2} and Proposition~\ref{pexsa5.3}. 
This  then yields  the explicit sufficiency condition on the degeneracy parameter $\delta>2-(d-d_{\!H})/2$ for self-adjointness.
Although there is a large literature on the Hardy constant it appears to be confined to domains with smooth boundaries or satisfying convexity properties.
Little appears to be known about the value of the Hardy constant for domains with rough boundaries.
Nevertheless one can derive some general properties of $b_\delta$ which lead to a criterion for it to have the standard value 
$(d-d_{\!H}+\delta-2)/2$.
In addition one can demonstrate by explicit example  that this is not always the case. 
In fact the Hardy constant can be arbitrarily small (see Example~\ref{exsa6.1}). 

\smallskip

First we examine the general properties of $b_\delta$.
The most complicated one  to establish   in the current setting is that $b_\delta$ is bounded from above by the standard value.
Note that the following proposition applies equally well to convex domains or Lipschitz domains.
In both these cases it establishes that the boundary Hardy constant  satisfies the upper bound $b_\delta\leq (\delta-1)/2$ for $\delta>1$.
This result was already cited for convex domains in Example~\ref{exsa5.4}.
The general result does, however,  require  uniformity of the domain  if $d_{\!H}\in[\,d-1,d\,\rangle$.

\begin{prop}\label{psa6.3}
Assume $\Omega$ is a uniform  domain  with an Ahlfors $s$-regular boundary $\Gamma$.
Further assume that if  $\delta>2-(d-s)$ then there is  a $b_{\delta,r}>0$ such that 
\[
\|d_\Gamma^{\,\delta/2}\,\nabla\psi\|_2\geq b_{\delta,r}\,\|d_\Gamma^{\,\delta/2-1}\psi\|_2
\]
for all $\psi\in C_c^1(\Gamma_{\!\!r})$.
It follows, since $s=d_{\!H}$, that  the Hardy constant $b_\delta$ satisfies
\begin{equation}
b_\delta\leq (d-d_{\!H}+\delta-2)/2
\;.\label{esa6.10}
\end{equation}
\end{prop}
\proof\
It is evident that one may assume
\begin{equation}
b_{\delta,r}\leq b_\delta=\inf\Big\{\|d_\Gamma^{\,\delta/2}\,\nabla\psi\|_2/\|d_\Gamma^{\,\delta/2-1}\psi\|_2:\psi\in C_c^1(\Gamma_{\!\!r}),\,r>0\Big\}
\;.\label{esa6.2}
\end{equation}
Moreover, since $b_\delta$ is the supremum over all possible choices of $b_{\delta, r}$ for small $r>0$ it suffices to prove that $b_{\delta,r}$
satisfies the bound (\ref{esa6.10}).
But this can be established  by modification of a lemma of Adam Ward \cite{War} (see also \cite{War1}, Lemma~5.1).

\begin{lemma}\label{lsa6.1}
For each $\beta\geq0$
\begin{equation}
b_{\delta,r}\leq |(\beta+\delta-2)/2|+\|d_\Gamma^{\,1-\beta/2}\nabla\psi\|_2/\|d_\Gamma^{\,-\beta/2}\psi\|_2
\;.\label{esa6.3}
\end{equation}
for all $\psi\in C_c^1(\Gamma_{\!\!r})$.
\end{lemma}
\proof\
Set $\alpha=\beta+\delta-2$.
Then replace $\psi$ in (\ref{esa6.3}) by $d_\Gamma^{\,-\alpha/2}\psi$.
It follows from the Leibniz rule and the triangle inequality that 
\begin{eqnarray*}
\|d_\Gamma^{\,\delta/2}\nabla(d_\Gamma^{\,-\alpha/2}\psi)\|_2&\leq&\|d_\Gamma^{\,\delta/2}(\nabla d_\Gamma^{\,-\alpha/2})\psi\|_2+\|d_\Gamma^{\,\delta/2}d_\Gamma^{\,-\alpha/2}(\nabla \psi)\|_2\\[5pt]
&=&|(\beta+\delta-2)/2|\,\|d_\Gamma^{\,-\beta/2}\psi\|_2+\|d_\Gamma^{\,1-\beta/2}(\nabla\psi)\|_2
\;.
\end{eqnarray*}
Similarly $\|d_\Gamma^{\,\delta/2-1}d_\Gamma^{\,-\alpha/2}\psi\|_2=\|d_\Gamma^{\,-\beta/2}\psi\|_2$.
The desired conclusion follows immediately.\hfill$\Box$

\bigskip

The next step in the proof of Proposition~\ref{psa6.3}  is to construct a sequence of $\psi_n\in C_c^1(\Gamma_{\!\!r})$ such that if $\beta=d-d_{\!H}$ then the numerator in the last term of (\ref{esa6.2})
 is bounded uniformly in $n$ but the denominator diverges as $n\to\infty$.
 Once this is achieved one immediately deduces the bound (\ref{esa6.10}).
 The construction is based on some estimates  which are a consequence of the uniformity of $\Omega$ and the Ahlfors regularity of $\Gamma$.
 
 Let $A=\Gamma\cap B(x_0;R)$, with $x_0\in \Gamma$ and $R>0$ and set $A_r=\{x\in\overline \Omega:d_A(x)<r\}$ where $d_A$ is the Euclidean distance to the set $A$.
 Then there are $\kappa,\kappa'>0$ such that 
 \begin{equation}
 \kappa'r^{(d-d_{\!H})}\leq |A_r|\leq \kappa\, r^{(d-d_{\!H})}
 \label{esa6.30}
 \end{equation}
 for all small $r>0$.
 These estimates are established in Section~2 of \cite{LR}.
 Their proof is based on ideas of Salli \cite{Sall}.
 The upper bound only requires the Ahlfors regularity of $\Gamma$ but if $s\in [\,d-1,d\,\rangle$ the lower bound also requires  the  uniformity property.

 The $\psi_n$ are now defined with the aid of the $W^{1,\infty}(0,\infty)$-functions $\xi_n$ used in the proof of Corollary~\ref{csa2.1}.
On $[0,1]$ one has $\xi_n(t)=0$ if $t<1/n$, $\xi_n(t)=1$ if $t>1$ and $\xi_n(t)=\log(nt)/\log n$ if $1/n\leq t\leq 1$.
Thus $\xi_n(1)=1$ and $\xi_n(t)\geq 1/2$ if $t\in [n^{-1/2},1]$.
Then fix a decreasing function  $\chi\in C^1(0,1)$ with $\chi(t)=1$ if $t\in[0,1/2]$, $\chi(1)=0$ and $|\chi'|\leq 4$.
Finally define $\psi_n=(\xi_n\circ(r^{-1}d_\Gamma))(\chi\circ(r^{-1}d_A ))$.

It follows from this construction that  $0\leq\psi_n\leq 1$, $\supp\psi_n\subset A_r\subset\Gamma_{\!\!r}$ and  the $\psi_n$ converge pointwise to $\chi\circ(r^{-1}d_A)$ as $n\to\infty$.
Further $(\xi_n\circ(r^{-1}d_\Gamma))\geq 1/2$ if $d_\Gamma>rn^2$ and $(\chi\circ(r^{-1}d_A ))\geq1$ if $d_A<r/2$.
In addition $d_A\geq d_\Gamma$.
Therefore, setting $\beta=d-d_{\!H}$, 
\[
\|d_\Gamma^{\,-\beta/2}\psi_n\|_2=\int d_\Gamma^{\,-\beta}\,|\psi_n|^2\geq (1/4)\int d_A^{\,-\beta}\,\one_{D_{r,n}}
\]
where $D_{r,n}=\{x\in\Omega:r/n^{1/2}\leq d_A(x)\leq r/2\,,\,r/n^{1/2}\leq d_\Gamma(x)\leq r\}$.
Then since $d_A^{\,-\beta}=(r/2)^{-\beta}(1+\beta\int^1_{d_A(x)/r}t^{-(\beta+1)})$ one has 
\begin{eqnarray*}
\int d_\Gamma^{\,-\beta}\,|\psi_n|^2&\geq &(1/4)(r/2)^{-\beta}\int \one_{D_{r,n}}\Big(1+\beta\int^1_{d_A(x)/r}t^{-(\beta+1)}\Big)\\[5pt]
&\geq &(1/4)(r/2)^{-\beta}|D_{r,n}| +(1/4)\,\beta\int^1_{r n^{-1/2}}dt\,t^{-1}((rt/2)^{-\beta}|D_{r,t,n}|)
\end{eqnarray*}
where $D_{r,t, n}=\{x\in\Omega:r/n^{1/2}\leq d_A(x)\leq tr/2\,,\,r/n^{1/2}\leq d_\Gamma(x)\leq r\}$.
But it follows from the volume estimates (\ref{esa6.30}) that $\lim_{n\to\infty}((r/2)^{-\beta}|D_{r,n}|)=((r/2)^{-\beta}|A_{r/2}|\geq \kappa'$ and 
$\lim_{n\to\infty}((rt/2)^{-\beta}|D_{r,t,n}|)=((rt/2)^{-\beta}|A_{rt/2}|)\geq \kappa'$
uniformly for $t\in\langle0,1\rangle$.
Since the integral of $t^{-1}$ is divergent at the origin one concludes that $\lim_{n\to\infty}\|d_\Gamma^{\,-\beta/2}\psi_n\|_2=\infty$.
This is the first step in handling the last term in (\ref{esa6.3}).

The second step is to estimate $\|d_\Gamma^{\,1-\beta/2}(\nabla\psi_n)\|_2^2=\int d_\Gamma^{\,2-\beta}|\nabla\psi_n|^2$.
First observe that
\[
|\nabla\psi_n|^2\leq 2r^{-2}|(\xi'_n\circ(r^{-1}d_\Gamma))|^2|(\chi\circ(r^{-1}d_A ))|^2
+2r^{-2}|(\xi_n\circ(r^{-1}d_\Gamma))|^2|(\chi'\circ(r^{-1}d_A ))|^2
\]
by the Leibniz rule and the Cauchy-Schwarz inequality.
Denote the integrals involving the first and second terms on the right hand side by $I_1$ and $I_2$, respectively.
Then since $|\xi_n|\leq 1$ one has 
\begin{eqnarray*}
I_2&\leq &8r^{-2}\int d_\Gamma^{\,2-\beta} |(\xi_n\circ(r^{-1}d_\Gamma))|^2\,\one_{\{x:r/2\leq d_A(x)\leq r\}}\\[5pt]
&\leq&8r^{-2}\int d_\Gamma^{\,2-\beta}\,\one_{\{x:r/2\leq d_A(x)\leq r\}}
\leq 8(r^{-\beta}|A_r|)\leq 8\kappa
\end{eqnarray*}
if $\beta\leq 2$.
Alternatively, if $\beta>2$ one deduces that 
\[
I_2\leq 8r^{-2}\int d_\Gamma^{\,2-\beta}\, \one_{D'_{r,n}}\leq 8\Big(r^{-\beta}\,|D'_{r,n}|+(\beta-2)\int^1_{n^{-1}}dt\,t\,((rt)^{-\beta}\,|D'_{r,t,n}|\Big)
\]
where $D'_{r,t,n}=\{x\in\Omega:rn^{-1}\leq d_\Gamma(x)\leq rt\,,\, r/2\leq d_A(x)\leq r\}$ and $D'_{r,n}=D'_{r,1,n}$.
As before one has $\lim_{n\to\infty}(r^{-\beta}|D'_{r,n}|)=(r^{-\beta}|A_{r}|)\leq \kappa$ but the estimate on $(rt)^{-\beta}|D'_{r,t,n}|$ is more delicate.
First note that $D'_{r,t,n}\subset\{x\in\Omega:0\leq d_\Gamma(x)\leq rt\,,\,r/2\leq d_A(x)\leq r\}$.
But if $t<1/2$ then the condition $d_\Gamma(x)\leq rt$ eliminates all those $x\in A_r$ such that $d_A(x)=d_\Gamma(x)$.
Therefore $D'_{r,t,n}\subset {\hat A}_{rt}$ where $\hat A$ is a bounded subset of $\Gamma$ determined by $A$.
This can be explicitly verified as follows.

Since $A=\Gamma\cap B(x_0\,;R)$ one clearly has  $D'_{r,t,n}\subset \Omega\cap B(x_0\,;R+r)$.
More specifically
\begin{eqnarray*}
D'_{r,t,n}&\subseteq &( \Omega\cap B(x_0\,;R+r))\cap\{x\in\Omega: d_A(x)\geq r/2\,,\, d_\Gamma(x)\leq rt\}\\[5pt]
&\subseteq &( \Omega\cap B(x_0\,;R+r)\backslash B(x_0\,;R))\cap\{x\in\Omega:  d_\Gamma(x)\leq rt\}
\end{eqnarray*}
for all small $t$.
The second inclusion follows since for $t<1/2$ the combination of the condition $d_A(x)\geq r/2$ and $d_\Gamma(x)<rt$ implies that $x$ is in the complement of $B(x_0\,;R)$.
Now it follows that one can choose $\hat A=\Gamma\cap B(x_0\,;R+r)$.
Therefore $(rt)^{-\beta}|D'_{r,t,n}|\leq (rt)^{-\beta}|{\hat A}_{rt}|$ is bounded uniformly for all $n\geq 1$ and $t\leq 1$ by (\ref{esa6.30}).
Consequently $I_2$ is uniformly bounded in $n$.

Next we have to estimate the integral $I_1$.
But 
\[
I_1=2r^{-2}\int d_\Gamma^{\,2-\beta}|(\xi'_n\circ(r^{-1}d_\Gamma))|^2|(\chi\circ(r^{-1}d_A ))|^2\leq 2(\log n)^{-2}\int d_\Gamma^{\,-\beta}\,\one_{D''_{r,n}}
\]
where one now has  $D''_{r,n}=\{x\in\Omega:rn^{-1}\leq d_\Gamma(x)\leq r\,,\,0< d_A(x)\leq r\}$.
Arguing as above
\[
I_1\leq 2(\log n)^{-2}\Big(|D''_{r,n}|+\beta\int^1_{n^{-1}} dt\,t^{-1}((rt)^{-\beta}\,|D''_{r,t,n}|\Big)
\]
with $D''_{r,t,n}=\{x\in\Omega:rn^{-1}\leq d_\Gamma(x)\leq rt\,,\, 0< d_A(x)\leq r\}$.
The the volume estimates  (\ref{esa6.30}) establish that $\sup_{n\geq1}(r^{-\beta}|D''_{r,n}|)\leq \kappa$ and  $\sup_{n\geq1}((rt)^{-\beta}|D''_{r,t,n}|)\leq \kappa$.
But the integral gives a factor $\log n$ and one obtains a bound $I_1\leq b\, (\log n)^{-1}$ with $b>0$.
Hence one concludes that
\[
 \|d_\Gamma^{\,1-\beta/2}(\nabla\psi_n)\|_2=\int d_\Gamma^{\,2-\beta}|\nabla\psi_n|^2\leq a+b\,(\log n)^{-1}
 \]
 with $a,b>0$.
 This completes the proof of  Proposition~\ref{psa6.3}.
 \hfill$\Box$
 
 \bigskip

It follows from the upper bound of  Proposition~\ref{psa6.3} combined with the lower bound of  Proposition~2.9 of \cite{Rob15}  that for the domains covered by 
Corollary~\ref{cexsa5.2} and Proposition~\ref{pexsa5.3}
the optimal Hardy constant always has the standard value $b_\delta=(d-d_{\!H}+\delta-2)/2$.
In the current setting this is not necessarily the situation.
Nonetheless one can draw some interesting conclusions from 
general properties of  $b_\delta$.

The next proposition collects three related properties which are   independent of the particular characteristics of $\Omega$ and $\Gamma$.

\begin{prop}\label{lsa6.2}
Assume $\Omega$ is such that the weighted Hardy inequality 
\[
\|d_\Gamma^{\,\delta/2}\,\nabla\psi\|_2\geq b_{\delta,r}\,\|d_\Gamma^{\,\delta/2-1}\psi\|_2
\]
is valid for all $\psi\in C_c^1(\Gamma_{\!\!r})$ and all $\delta\in I=\langle 2-(d-d_{\!H}), 2\,]$.

\smallskip

It follows that 
\begin{tabel}
\item\label{lsa6.2-1}
$\delta\in I\mapsto b_\delta+\delta/2$ is non-decreasing and $b_2\geq b_\delta-(2-\delta)/2$ for all $\delta\in I$,
\item\label{lsa6.2-2}
if $b_2\geq (d-d_{\!H})/2$ then $|b_2-b_\delta|\leq (2-\delta)/2$ for all $\delta\in I$,
\item\label{lsa6.2-3}
if, in addition,  $b_\delta\leq (d-d_{\!H}+\delta-2)/2$  then $b_2=b_\delta-(2-\delta)/2$ for all $\delta\in I$.
\end{tabel}
\end{prop}
\proof\
{\ref{lsa6.2-1}.}
Fix $\delta_1,\delta_2\in I$  with $\delta_1<\delta_2$.
Then
\begin{eqnarray*}
\|d_\Gamma^{\,\delta_2/2-1}\psi\|_2&=&\|d_\Gamma^{\,\delta_1/2-1}(d_\Gamma^{\,(\delta_2-\delta_1)/2}\psi)\|_2\\[5pt]
&\leq& b_{\delta_1,r}^{\,-1}\,\|d_\Gamma^{\,\delta_1/2}\,\nabla (d_\Gamma^{\,(\delta_2-\delta_1)/2}\psi)\|_2\\[5pt]
&\leq& b_{\delta_1,r}^{\,-1}\Big(\|d_\Gamma^{\,\delta_2/2}(\nabla \psi)\|_2+((\delta_2-\delta_1)/2)\,\|d_\Gamma^{\,\delta_2/2-1}\psi\|_2\Big)\\[5pt]
&\leq &b_{\delta_1,r}^{\,-1}\Big(1+((\delta_2-\delta_1)/2)\,b_{\delta_2,r}^{-1}\Big)\,\|d_\Gamma^{\,\delta_2/2}(\nabla \psi)\|_2
\end{eqnarray*}
for all $\psi\in C_c^1(\Gamma_{\!\!r})$.
Therefore 
\[
b_{\delta_1,r}\Big(1+((\delta_2-\delta_1)/2)\,b_{\delta_2,r}^{-1}\Big)^{-1}\leq b_{\delta_2,r}
\]
or, equivalently,
$b_{\delta_1,r}+\delta_1/2\leq b_{\delta_2,r}+\delta_2/2$.
Thus taking the supremum over the choice of the $b_{\delta_j,r}$ one has $b_{\delta_1}+\delta_1/2\leq b_{\delta_2}+\delta_2/2$.
Hence  $\delta\in I\mapsto b_\delta+\delta/2$ is  non-decreasing.

\smallskip

\noindent{\ref{lsa6.2-2}.}
First it follows from  I  that $b_2+1\geq b_\delta+\delta/2$ for all $\delta\in I$.
Therefore $b_2-b_\delta\geq -(2-\delta)/2$.

Secondly, $(2-\delta)/2<(d-d_{\!H})/2$ since $\delta\in I$.
Hence  $b_2-(2-\delta)/2>b_2-(d-d_{\!H})/2\geq 0$
by assumption.
Thirdly, it follows from the weighted  Hardy inequality  that
\[
 b_{2,r}\,\|d_\Gamma^{\,\delta/2-1}\psi\|_2\leq \|d_\Gamma\nabla(d_\Gamma^{\,\delta/2-1}\psi)\|_2\leq (1-\delta/2)\|d_\Gamma\,d_\Gamma^{\,\delta/2-2}\psi\|_2+\|d_\Gamma^{\,\delta/2}\,\nabla\psi\|_2
 \;.
 \]
But by the preceding  one may choose $r$ sufficiently small that  $(2-\delta)/2<b_{2,r}$.
Hence
\[
(b_{2,r}-(2-\delta)/2)\,\|d_\Gamma^{\,\delta/2-1}\psi\|_2\leq \|d_\Gamma^{\,\delta/2}\,\nabla\psi\|_2
\]
and  $b_\delta\geq b_2-(2-\delta)/2$.
Therefore $b_2-b_\delta\leq (2-\delta)/2$.
Thus $|b_2-b_\delta|\leq (2-\delta)/2$ for all $\delta\in I$.
 \smallskip

\noindent{\ref{lsa6.2-3}.}
If $b_2\geq (d-d_{\!H})/2$ and $b_\delta\leq (d-d_{\!H}+\delta-2)/2$ then it immediately follows that  $b_\delta\leq (d-d_{\!H})/2-(2-\delta)/2\leq b_2-(2-\delta)/2$.
Therefore $b_\delta\leq b_2-(2-\delta)/2$.
But the converse  was established in  the proof of II.
Hence one must have an equality.
\hfill$\Box$

\bigskip

The final assumption $b_\delta\leq (d-d_{\!H}+\delta-2)/2$ in Proposition~\ref{lsa6.2} was established in Proposition~\ref{psa6.3}.
Therefore under the assumptions of that proposition one concludes that the bound $b_2\geq (d-d_{\!H})/2$ implies that $b_\delta=b_2-(2-\delta)/2$ and 
$b_\delta $ is a strictly increasing function of $\delta$ on the interval $\langle2-(d-d_{\!H}),2\,]$.
This observation is used in the proof of the following 
more precise  version of Theorem~\ref{tsa5.1}.
In contrast to that theorem we do not have to assume explicitly  the validity of the weighted Hardy inequality as it is a consequence of Lehrb{\"a}ck's theorem,
Proposition~\ref{psa6.2}.
  
  \begin{thm}\label{tsa6.1}
Assume that $\Omega$ is a uniform domain whose boundary $\Gamma$ is Ahlfors $s$-regular.
Further assume  that  the coefficients $C$ of $H$ satisfy the boundary condition $( \ref{esa5.1})$.
Then one has the following:
\begin{tabel}
\item \label{tsa6.1-1}
 there is a unique $\delta_c\in[\, 2-(d-d_{\!H})/2,2\rangle$ such that $H$ is self-adjoint for all $\delta> \delta_c$,
 \item \label{tsa6.1-2}
 $b_\delta=(d-d_{\!H}+\delta-2)/2$ for all  $\delta\in\langle  2-(d-d_{\!H})/2,2\,]$ if and only if  $b_2\geq(d-d_{\!H})/2$, 
 and if these conditions are satisfied then $b_2=(d-d_{\!H})/2$ and  $\delta_c=2-(d-d_{\!H})/2$.
 \end{tabel}
 \end{thm}
 \proof\
 \ref{tsa6.1-1}. First it  follows from Corollary~\ref{csa2.3} that $H$ is self-adjoint if $\delta\geq2$.

Secondly, the weighted Hardy inequality  is valid on a boundary layer $\Gamma_{\!\!r}$  for all  $\delta>2-(d-d_{\!H})$ by Proposition~\ref{psa6.2}.
Therefore it follows from
Theorem~\ref{tsa5.1} that  the condition $(2-\delta)/2<b_\delta$  is sufficient for the self-adjointness of $H$ for 
  $\delta\in \langle 2-(d-d_{\!H}), 2\,]$.
  Setting $B_\delta=b_\delta+\delta/2$ this sufficiency condition becomes $B_\delta>1$.
But  $b_\delta>0$ for $\delta>2-(d-d_{\!H})$ by Proposition~\ref{psa6.2}.
Therefore $B_\delta>\delta/2$.
 In particular $B_2>1$.
 
 Thirdly, $b_\delta\leq (d-d_{\!H}+\delta-2)/2$ by Proposition~\ref{psa6.3}.
 Therefore  $B_\delta\leq (d-d_{\!H}+2\delta-2)/2$. 
 In particular $B_{2-(d-d_{\!H})/2}\leq 1$.
 
 Finally $b_\delta=b_2-(2-\delta)/2$ for   $\delta\in \langle 2-(d-d_{\!H}),2\,]$ by the remark preceding the theorem.
Therefore  there  must be a unique  $\delta_c$ in this range for which $B_{\delta_c}=1$.
Thus if  $\delta>\delta_c$ then $H$ is self-adjoint by Theorem~\ref{tsa5.1}.

\smallskip

\noindent{ \ref{tsa6.1-2}.} 
Assume  $b_2\geq (d-d_{\!H})/2$.
It follows from Proposition~\ref{psa6.3} that $b_\delta\leq (d-d_{\!H}+\delta-2)/2$.
In particular $b_2\leq (d-d_{\!H})/2$.
Therefore $b_2= (d-d_{\!H})/2$.
But $b_\delta=b_2-(2-\delta)/2$ by the remark preceding the theorem.
Hence $b_\delta=(d-d_{\!H}+\delta-2)/2$.
Conversely if  $b_\delta=(d-d_{\!H}+\delta-2)/2$ then $b_2=(d-d_{\!H})/2$.

Finally if these equivalent conditions are satisfied then  $b_\delta>(2-\delta)/2$ is equivalent to the condition $\delta>2-(d-d_{\!H})/2$.
Therefore one has  $\delta_c=2-(d-d_{\!H})/2$.
\hfill$\Box$

\bigskip

It now follows from the remark preceding the theorem and 
Statement~\ref{tsa6.1-2} of the theorem 
that one has   the following reformulation of the self-adjointness result for the standard case.

\begin{cor}\label{csa6.11}
Assume $\Omega$ is a uniform domain with an Ahlfors $s$-regular boundary $\Gamma$
and that  the coefficients $C$ of $H$ satisfy the boundary condition $( \ref{esa5.1})$.

\smallskip
If $b_2\geq (d-d_{\!H})/2$ then $H$ is self-adjoint for all $\delta>2-(d-d_{\!H})/2$
and and the Hardy constant $b_\delta=(d-d_{\!H}+\delta-2)/2$ for all  $\delta\in\langle  2-(d-d_{\!H})/2,2\,]$.
\end{cor}

Next consider the non-standard case under the ongoing uniformity and Ahlfors regularity assumptions.
Since $b_\delta\leq (d-d_{\!H}+\delta-2)/2$ for $\delta>2-(d-d_{\!H})$ by Proposition~\ref{lsa6.2} the non-standard case corresponds to the condition 
$b_\delta<(d-d_{\!H}+\delta-2)/2$.
The following example demonstrates that this can occur and one can even have $b_\delta$ arbitrarily small.

\begin{exam}\label{exsa6.1}
First let $B_R$ denote a ball of radius $R$.
 Secondly let $\psi\in C_c^1(B_R)$ denote a function which is one on  the concentric ball $B_{(1-2\varepsilon)R}$, zero on $B_{(1-\varepsilon)R}$ and 
such that  $\|\nabla\psi\|_\infty\leq 2/(\varepsilon R)$ where  $\varepsilon\in \langle0,1/4\rangle$.
 Then
 \[
 \int_{B_R}d_\Gamma^{\,\delta-2}|\psi|^2\geq \int_{B_{R/2}}d_\Gamma^{\,\delta-2}|\psi|^2\geq aR^{\,d+\delta-2}
 \;\;\;\;{\rm and}\;\;\;\;\int_{B_R}d_\Gamma^{\,\delta}\,|\nabla\psi|^2\leq bR^{\,d+\delta-2} \varepsilon^{\delta-1}
 \]
 with $a,b>0$ where the values of $a$ and $b$ are independent of $R,\varepsilon$ and $\delta$.
 Therefore
 \[
 \int_{B_R}d_\Gamma^{\,\delta}\,|\nabla\psi|^2\Big/\int_{B_R}d_\Gamma^{\,\delta-2}|\psi|^2\leq (b/a)\,\varepsilon^{\delta-1}
 \]
 and if $\delta>1$ this ratio tends to zero as $\varepsilon\to0$.

Secondly, let  $B_k$ denote  the concentric ball with radius $2^{-k}R$ where $k=0,1,\ldots$.
Thus $B_0=B_R$. 
Then let  $\psi_k\in C_c^1(B_k)$ be the scaled function with $\psi_k(x)=\psi(2^{k}x)$ for $x\in B_0$.
It follow immediately by scaling invariance that the above ratio is unchanged by the replacement $B_0\to B_k$ and $\psi\to\psi_k$.
\end{exam}

\vspace{-4mm}
\begin{wrapfigure}{r}{5cm}
\begin{center}
\includegraphics[width=4cm,keepaspectratio]{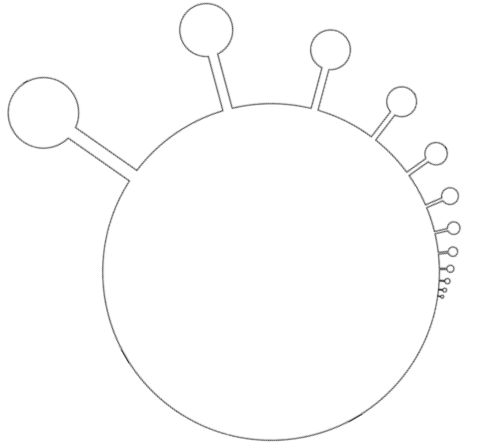}
\end{center}
\end{wrapfigure}

Thirdly, let $B$ denote a large ball and attach (a family of translates of)  the  balls $B_k$ to $B$  by narrow tunnels of length
$2^{-k}R$ and width $\varepsilon\,(2^{-k}R)$.
The balls and tunnels are understood to lie in the exterior of $B$ with  the attachments   separated from each other such that no pair overlap, as indicated in the illustration.
Then let $\Omega$ denote the open interior of the `decorated' ball.
It follows that $\Omega$ is a uniform domain with an Ahlfors $(d-1)$-regular boundary.
In particular $d_{\!H}=d-1$.
Hence  the standard value $(d-d_{\!H}+\delta-2)/2$ of the Hardy constant would be $(\delta-1)/2>0$.
 
 Finally let $\Psi_{\!n}$ be the function with $\supp\Psi_{\!n}=\bigcup_{k\geq n}B_k$ such that $\Psi_{\!n}|_{B_k}=\psi_k$.
 It then follows that for each small $r>0$ there is an $n$ such that $\supp\Psi_{\!n}\subset \Gamma_{\!\!r}$.
 Therefore $b_{\delta,r}\leq 2 (b/a)\varepsilon^{\delta-1}$ where the factor~$2$ is added to compensate any small adjustment needed to account for the effect of the tunnels.
 Hence $b_\delta\leq 2(b/a)\varepsilon^{\delta-1}<(\delta-1)/2$ for all  sufficiently small $\varepsilon$.
 Therefore $b_\delta$ does not attain the standard value and can be made arbitrarily small by appropriate choice of~$\varepsilon$.
 
 \bigskip
 
  Although this example is rather special it does indicate a general property.
 The value of the parameter $\sigma$ which enters the definition of uniform domains is governed by the ratio $2/\varepsilon$ of the length  of the tunnels
 and the width of the tunnels.
 In particular as $\varepsilon$ decreases the uniformity parameter increases.
 In addition the Ahlfors parameter $a$ governing the regularity of the boundary has a dependence on the width of the tunnels.
Therefore  one must expect that the Hardy constant   depends both  on  the degree of uniformity
 and also  the Ahlfors regularity parameter.

 \smallskip
 
 Finally we note that the value $b_2$ of the Hardy constant has a particular significance.
 First since $b_\delta\leq (d-d_{\!H}+\delta-2)/2$ if $\delta>2-(d-d_{\!H})$ by Proposition~\ref{psa6.3} one must have $b_2\leq (d-d_{\!H})/2$.
 Moreover, $b_\delta=(d-d_{\!H}+\delta-2)/2$ if and only if $b_2=(d-d_{\!H})/2$.
 Secondly, the weighted Hardy inequality with $\delta=2$ states that
 \[
 (\psi, H_2\psi)\geq b_{2,r}^{\,2}(\psi, \psi)
 \]
 for all $\psi\in C_c^2(\Gamma_{\!\!r})$ where $H_2\psi=-\divv(d_\Gamma^{\,2}\,\nabla\psi)$.
 Thus
 \[
 b_{2,r}^{\,2}\leq \inf\{(\psi, H_2\psi):\psi\in C_c^2(\Gamma_{\!\!r})\,,\; \|\psi\|_2=1\}
 \;.
 \]
 Therefore the supremum of the possible $ b_{2,r}^{\,2}$ is equal to the infimum of the spectrum of $H_2$ acting on $L_2(\Gamma_{\!\!r})$
 and  $b_2^{\,2}$ is the infimum of the spectra of $H_2$ acting on $L_2(\Gamma_{\!\!r})$ for all small $r$.
 In particular the bottom of the spectrum is $(d-d_{\!H})^2/4$ if and only if $b_\delta$ has the standard value.
 In the non-standard case the bottom of the spectrum is strictly smaller than $(d-d_{\!H})^2/4$  and as the example shows it can be arbitrarily small.
This observation might  provide a practical method in special cases, such as convex domains, of confirming that the Hardy constant has the standard value.

\section{Rellich inequalities}

Theorems~\ref{tsa5.1} and \ref{tsa6.1} demonstrate that the weighted Hardy inequality is the essential ingredient
in the derivation of self-adjointness of the degenerate elliptic operator $H$.
In contrast  in the earlier paper \cite{Rob15}  a weighted Rellich inequality for functions supported near the boundary
was the crucial element.
Although the Rellich inequality is no longer needed for the self-adjointness problem it is of interest that it can nevertheless
 be derived  in the broader setting of John-Ahlfors domains.

First we derive a  Hardy inequality for the quadratic form $h$ associated with $H$.
Note that Lehrb{\"a}ck's result, Proposition~\ref{psa6.2}, establishes the weighted Hardy inequality (\ref{esa5.3}) on a boundary layer
for John domains with Ahlfors regular boundaries so the following statements follow in this setting.

 \begin{lemma}\label{lsa7.1}
 Assume that  the coefficients $C$ of $H$ satisfy the boundary condition $( \ref{esa5.1})$ and that the weighted 
 Hardy inequality $(\ref{esa5.3})$ is valid on the boundary layer $\Gamma_{\!\!r}$ with $\delta> 2-(d-d_{\!H})$.

\smallskip

Then there are $s\in\langle0,r\rangle$  and $a_{\delta, s}\in\langle0,b_{\delta,r}\rangle$ such that the  Hardy inequality
\[
h(\psi)\geq a_{\delta,s}^{\,2}\,\|c^{1/2}d_\Gamma^{\,\delta/2-1}\,\psi\|_2^2
\]
is valid for all $\psi\in D(h)$ with $\supp \psi\subset \Gamma_{\!\!r}$.
$\psi\in C_c^1(\Gamma_{\!\!s})$. 
Moreover $a_\delta$,  the supremum of the possible $a_{\delta,s}$,  is equal to $b_\delta$,
the Hardy constant of  $(\ref{esa5.3})$. 
 \end{lemma}
\proof\ $\;$
First it follows from the boundary condition (\ref{esa5.2}) that 
\[
h(\psi)\geq \sigma_{\!r}\|c\,d_\Gamma^{\,\delta/2}(\nabla\psi)\|_2^2
\]
for all $\psi\in C_c^1(\Gamma_{\!\!r})$.
Hence \begin{eqnarray*}
h(\psi)^{1/2}&\geq& \sigma_{\!r}^{1/2}\Big(\|d_\Gamma^{\,\delta/2}\,\nabla(c^{1/2}\psi)\|_2-a_c\,\|c^{1/2}d_\Gamma^{\,\delta/2}\,\psi\|_2\Big)\\[5pt]
&\geq& \sigma_{\!r}^{1/2}\Big(b_{\delta,s}\|c^{1/2}d_\Gamma^{\,\delta/2-1}\psi\|_2-sa_c\,\|c^{1/2}d_\Gamma^{\,\delta/2-1}\,\psi\|_2\Big)
=(\sigma_{\!r}^{1/2}\,a_{\delta,s})\,\|c^{1/2}d_\Gamma^{\,\delta/2-1}\psi\|_2
\end{eqnarray*}
where $a_{\delta,s}=b_{\delta,s}-sa_c$ and  $a_c=2^{-1}\|(\nabla c)/c\|_\infty$ with $\|\cdot\|_\infty$  the $L_\infty(\Gamma_{\!\!r})$-norm.
Note that the  $b_{\delta,s}>0$ may be chosen such  that $b_{\delta,s}$ converges upward to $b_\delta$ as $s\to0$.
Therefore one may assume $a_{\delta,s}>0$ for all small $s>0$
Then  $a_{\delta,s}\to b_\delta$ as $s\to0$ and since
$\sigma_{\!r}\to1$ as $r\to0$ the proof is complete for $\psi\in  C_c^1(\Gamma_{\!\!r})$.
But the Hardy inequality extends by continuity to the $\psi\in D(h)$ with support in the boundary layer.
\hfill$\Box$

\bigskip

The Hardy inequality of the lemma now implies   the Rellich inequality  by  adaptation of the simpler part of  the reasoning in Section~2 of \cite{Rob12}.

 \begin{prop}\label{psa7.2}
 Assume that  the coefficients $C$ of $H$ satisfy the boundary condition $( \ref{esa5.1})$ and that the weighted 
 Hardy inequality $(\ref{esa5.3})$ is valid on $\Gamma_{\!\!r}$ with $\delta> 2-(d-d_{\!H})$.

\smallskip

If $b_\delta>(2-\delta)/2$ then there are $r>0$ and a map $s\in\langle0,r\rangle\mapsto B_{\delta,s}>0$ such that  the Rellich inequality
\begin{equation}
\|H\psi\|_2^2\geq B_{\delta,s}^{\,2}\,\|c\,d_\Gamma^{\,\delta-2}\psi\|_2^2
\label{ensa1}
\end{equation}
is valid for all $\psi$ in the domain $ D(H)$ of the  self-adjoint operator $H$ with $\supp\psi\subset\Gamma_{\!\!s}$. 
Moreover, the Rellich constant  $B_\delta$,  the supremum of the possible $B_{\delta,s}$, is given by
 $B_\delta=b_\delta^{\,2}-((2-\delta)/2)^2$.
\end{prop}
\proof\
First it follows from Lemma~\ref{lsa7.1} that there is an $s>0$  such that the Hardy inequality $h(\psi)\geq \|\chi_s\psi\|_2^2$
is valid for all $\psi\in C_c^\infty(\Gamma_{\!\!s})$ with $\chi_s=a_{\delta,s}(c^{1/2}d_\Gamma^{\,\delta/2-1})$.

Secondly, it  follows from the boundary condition (\ref{esa5.2}) that 
\[
\Gamma_{\!\!c}(\chi_s)\leq \tau_s\,(c\,d_\Gamma^{\,\delta})\,|\nabla\chi_s|^2
\;.
\]
Then a straightforward calculation establishes the eikonal inequality 
\begin{equation}
\Gamma_{\!\!c}(\chi_s)\leq {\gamma}_s\,\chi_s^4
\label{ensa2}
\end{equation}
where
${\gamma}_s=\tau_s\,((2-{\delta}_s)/(2\,a_{\delta,s}))^2$
with ${\delta}_s=\delta-sa_c$ and  we have used the notation of the proof of Lemma~\ref{lsa7.1}
(see  \cite{Agm1}, Theorem~1.4(ii)).

Thirdly, it follows from the basic identity (\ref{esa2.1}) that 
\[
(H\psi, \chi_s^2\psi)=h(\chi_s\psi)-(\psi,\Gamma_{\!\!c}(\chi_s)\psi)
\]
for all $\psi\in C_c^2(\Gamma_{\!\!s})$.
Therefore
\[
\|H\psi\|_2\,\|\chi_s^2\psi\|_2\geq (1-{\gamma}_s)\,\|\chi_s^2\psi\|_2^2
\]
for all $\psi\in C_c^2(\Gamma_{\!\!s})$ by the Hardy inequality and (\ref{ensa2}).
It follows immediately that if ${\gamma}_s<1$ then
\[
\|H\psi\|_2^2\geq (1-{\gamma}_s)^2\,\|\chi_s^2\psi\|_2^2=
a_{\delta,s}^{\,4}\,(1-{\gamma}_s)^2\,\|c\,d_\Gamma^{\,\delta-2}\psi\|_2^2
\;.
\]
After substituting the value of ${\chi}_s$ and rearranging one obtains the Rellich inequality (\ref{ensa1}) with
$B_{\delta,s}=a_{\delta,s}^{\,2}\,(1-{\gamma}_s)=a_{\delta,s}^{\,2}-\tau_s\,((2-{\delta}_s)/2)^2>0$ for all $\psi\in C_c^2(\Gamma_{\!\!r})$.

Now  $\tau_s$ and $2-\delta_s$ converge downward to $1$ and $2-\delta$, respectively, as $s\to0$.
Moreover $a_{\delta,s}$ converges upward to $b_\delta$ by Lemma~\ref{psa7.2}.
Therefore $B_{\delta,s}$ converges upward to  $ B_\delta=b_\delta^{\,2}-((2-{\delta})/2)^2$ as $s\to0$.
Since $B_{\delta,s}>0$ for all small $s>0$ it follows that $B_\delta>0$.
Hence $b_\delta>(2-\delta)/2$.
Conversely, if $b_\delta>(2-\delta)/2$ then $B_\delta>0$ and $B_{\delta, s}>0$  or, equivalently, $\gamma_s<1$ for all small $s>0$.

Finally if $b_\delta>(2-\delta)/2$ then $H$ is self-adjoint by Theorem~\ref{tsa5.1}.
Therefore $C_c^2(\Omega)$ is a core of $H$.
Then the Rellich inequality (\ref{ensa1}) extends from $C_c^2(\Gamma_{\!\!r})$ to all $\psi\in D(H)$  with support in
$\Gamma_{\!\!r}$.
%
This follows since $\psi$ can  be approximated by a sequence $\psi_n\in C_c^\infty(\Omega)$ and then the $\psi_n$ 
can be replaced by $\xi\psi_n$ where $\xi$ is a $C^\infty$-function with support in ${\overline\Gamma_{\!\!r}}$ which is equal to 
$1$ on $\Gamma_{\!\!s}$ with $s<r$.
The modified functions are  in $D(H)$ and  still approximate $\psi$ in the graph norm as a simple corollary of Proposition~2.1.II in \cite{Rob15}.
\hfill$\Box$

\section{Summary and comments}

The foregoing investigation developed from the earlier works \cite{RSi4} and \cite{LR} on Markov uniqueness.
The Markov property is equivalent to the parabolic diffusion equation
\[
\partial\varphi_t/\partial t+H\varphi_t=0
\]
having a unique weak solution on $L_1(\Omega)$ (see \cite{Dav14} or \cite{RSi4})
whilst self-adjointness of $H$ is equivalent to uniqueness on  $L_2(\Omega)$.
The main conclusion of \cite{LR} for uniform domains with Ahlfors regular boundaries
 was the equivalence of $L_1$-uniqueness with the condition $\delta\geq 2-(d-d_{\!H})$. 
In  particular $L_1$-uniqueness depends on the geometry of the domain only through the 
Hausdorff dimension of its boundary.
The  foregoing analysis indicates that the situation is more complicated for $L_2$-uniqueness.
Although we have only derived the  sufficiency condition $b_\delta>(2-\delta)/2$
it is plausible that $b_\delta\geq(2-\delta)/2$ is both necessary and sufficient for 
self-adjointness, i.e.\ for $L_2$-uniqueness.
There are not many examples for guidance.
One simple case is  $\Omega=\Ri^d\backslash \{0\}$  and the operator 
$H=-\sum^d_{k=1}\partial_k(|x|^\delta\partial_k)$.
If $\delta>2-d$ the Hardy inequality is valid with $b_\delta=(d+\delta-2)/2$ and the condition $b_\delta>(2-\delta)/2$ is equivalent
to $\delta>2-d/2$.
Since, however,  the operator is rotationally invariant it can be  established by taking radial coordinates and applying  the classical Weyl limit-point,
limit-circle, theory  that $H$ is essentially self-adjoint if and only if $\delta\geq 2-d/2$.
A closely related result was established in \cite{Rob15} for $C^2$-domains.
Theorem~\ref{tsa5.1} establishes that $\delta>3/2$ is sufficient for self-adjointness but Theorem~3.2 of \cite{Rob15} also established that
the condition $\delta\geq 3/2$ is necessary. 
The latter result did, however, require  some bounds on the derivatives of the coefficients of the operator.

One thing that is clear is that the Hardy constant $b_\delta$ corresponding to the weighted  Hardy inequality on boundary layers plays a significant role
in determining self-adjointness.
A similar conclusion was reached by Ward \cite{War} \cite{War1} in his analysis of Schr{\"o}dinger operators on domains
by extension of the arguments of Nenciu and Nenciu \cite{NeN1}.
In addition $b_\delta$  can have a quite complicated dependency on the geometry of the boundary.
For domains with $C^2$-boundaries the Hardy constant only depends on the Hausdorff dimension of the boundary.
It is given by  $b_\delta=(d-d_{\!H}+\delta-2)/2$. 
An expression  we have referred to as the standard value. 
This is also the case for the complement of lower dimensional $C^2$-domains or the complement of convex subsets.
 Example~\ref{exsa6.1} demonstrates, however,  that for uniform domains with Ahlfors regular boundaries 
$b_\delta$, and consequently the self-adjointness property,  can  depend on the regularity parameters governing the boundary.
 It appears unlikely that one could calculate $b_\delta$ exactly in such situations.
 Nevertheless   Propositions~\ref{psa6.3} and \ref{lsa6.2} do  give 
general properties which allow one to gain considerable information about the diffusion.
In particular $b_\delta$ is  bounded from above by the standard value and the
corresponding critical degeneracy  for self-adjointness is larger than the value $2-(d-d_{\!H})/2$ in the standard case.
This raises the question as to the minimal smoothness requirements on the domain and its boundary to ensure that $b_\delta$
attains the standard value. 
For example, is this the case for $C^{1,1}$-domains or, more generally, for Lipschitz domains.
What is the situation for convex domains?

Finally we have only considered symmetric diffusion operators but our arguments should extend to the non-symmetric operators with drift terms 
considered by Nenciu and Nenciu \cite{NeN}.
Following these authors the non-symmetric operators on $L_2(\Omega)$ can be reformulated as symmetric operators on weighted spaces $L_2(\Omega\,; \rho)$.
Then, however, one would need some control on the behaviour of the weights $\rho$ near the boundary as in Theorem~5.3 of \cite{NeN}.

\section*{Acknowledgements}
Over the last five years I have been fortunate to have had a friendly and informative correspondence with Juha Lehrb{\"a}ck  on many matters ranging
from mathematics to marathon running.
I am now indebted to him for confirming the validity of the key Proposition~\ref{psa6.2},  providing the illuminating Example~\ref{exsa6.1}
and commenting on various drafts of the current paper.

I would also like to thank Louisa Barnsley for the illustration of Example~\ref{exsa6.1} based on Juha's earlier sketch.

\end{document}